\RequirePackage{fixltx2e}
\documentclass[a4paper]{isp}

\usepackage[american]{babel}

\usepackage{graphicx}

\usepackage{paralist}

\usepackage[T1]{fontenc}

\usepackage{lmodern}
\usepackage{microtype}
\usepackage{amsmath}
\usepackage{booktabs}
\usepackage{url}
\usepackage{amssymb}
\usepackage{lineno}

\ifnum\pdfoutput>0
\usepackage[
bookmarks=false,
breaklinks=true,
colorlinks=true,
linkcolor=black,
citecolor=black,
urlcolor=black,
pdfpagelayout=SinglePage
]{hyperref}
\usepackage[all]{hypcap}
\else
\usepackage{hyperref}
\fi

\usepackage[capitalise,nameinlink]{cleveref}
\newtheorem{theorem}{Theorem}
\newtheorem{definition}{Definition}%[theorem]
\crefname{section}{Sect.}{Sect.}
\Crefname{section}{Section}{Sections}
\crefname{figure}{Fig.}{Fig.}
\Crefname{figure}{Figure}{Figures}
\Crefname{theorem}{Th.}{Th.}
\Crefname{theorem}{Theorem}{Theorems}
\usepackage{xspace}
\begin{document}
%\pdfgentounicode=1

%\mainmatter

\title{Implementation of a Wiener Chaos Expansion Method for the Numerical Solution of the Stochastic Generalized Kuramoto-Sivashinsky Equation driven by Brownian motion forcing
	}
%If Title is too long, use \titlerunning
%\titlerunning{Short Title}

%Single insitute
%Currently disabled
\iffalse
%\author{Sheldon Cooper \and Leonard Hofstadter \and Rajesh Koothrappali \and Howard Wolowitz}
\author{Victor Nijimbere}
%If there are too many authors, use \authorrunning
%\authorrunning{First Author et al.}
%\institute{ Division of Physics, Mathematics and Astronomy, California Institute of Technology, Pasadena, United States}
\institute{ Carleton University, Ottawa, Ontario, Canada}
\fi

%%Multiple insitutes
%%Multiple institutes are typeset as follows:
%\author{Sheldon Cooper\inst{1}\corauthor \and Leonard Hofstadter\inst{2} \and Rajesh Koothrappali\inst{1} \and Howard Wolowitz\inst{2} }
\author{Victor Nijimbere }
%
%%If there are too many authors, use \authorrunning
%%\authorrunning{First Author et al.}
%
%\institute{
% Division of Physics, Mathematics and Astronomy, California Institute of Technology, Pasadena, United States\\
%\and
%Division of Engineering and Applied Science, California Institute of Technology, Pasadena, United States\\
%\email{sheldon.cooper@caltech.edu}
%}
\institute{
Carleton University, Ottawa, Ontario, Canada}
			
\maketitle

%\linenumbers

\begin{abstract}
 Numerical computations based on the Wiener Chaos Expansion (WCE) are carried out to approximate the solutions of the stochastic generalized Kuramoto--Sivashinsky (SgKS) equation driven by Brownian motion forcing. In the assessment of the accuracy of the WCE based approximate numerical solutions, the WCE based solutions are contrasted with semi-analytical solutions, and the absolute and relative errors are evaluated. It is found that the absolute error is $O(\varsigma t)$, where $\varsigma$ is small constant and $t$ is the time variabe; and the relative error is order $10^{-2}$ or less. This demonstrates that numerical methods based on the WCE are powerful tools to solve the SgKS equation or other related stochastic evolution equations.
\end{abstract}

%%%%%%%%%%%%%%%%%%%%%%%%%%%%%%%%%%%%%%%%%%%%%%%%%%%%%%%%%%%%%%%%%%%%%%%%%%%%%%%
\section{Introduction}\label{sec:1}
%%%%%%%%%%%%%%%%%%%%%%%%%%%%%%%%%%%%%%%%%%%%%%%%%%%%%%%%%%%%%%%%%%%%%%%%%%%%%%%
The evolution of complex spatiotemporal mechanisms, in fluids for example, involving nonlinear interactions can be modeled by quasi-nonlinear stochatic partial differential equations (SPDEs). In the present study, we focus on the stochastic generalized Kuramoto-Sivashinsky (SgKS) equation driven by a time-dependent Brownian motion forcing
\begin{equation}
\frac{\partial u}{\partial t}=-u \frac{\partial u}{\partial x}-\kappa \frac{\partial^2 u}{\partial x^2}-\eta \frac{\partial^3 u}{\partial x^3}-\nu\frac{\partial^4 u}{\partial x^4}+\sigma \dot{W},
\label{eq:1}
\end{equation}
where $u$ depends on both time $t$ and position $x$, $\kappa$ is a positive parameter representing the relative importance of the effects due instabilities associated to energy production, $\eta$ is a parameter characterizing the effects due to wave dispersion, $\nu$ is the viscous damping coefficient associated to energy dissipation, the amplitude of the forcing $\sigma$ is a function of $x$,  $W$ is a scalar function of time $t$ denoting the Brownian motion and $\dot{W}$ its derivative with respect to time $t$ (Gaussian white noise). The stochastic forcing term, in general, represents irregular random scatters which may affect a system under consideration. In the absence of wave dispersion, $\eta=0$, the SgKS equation (\ref{eq:1}) becomes the stochastic Kuramoto-Sivashinsky (SKS) equation
\begin{equation}
\frac{\partial u}{\partial t}=-u \frac{\partial u}{\partial x}-\kappa \frac{\partial^2 u}{\partial x^2}-\nu\frac{\partial^4 u}{\partial x^4}+\sigma \dot{W}.
\label{eq:2}
\end{equation}

The SPDEs (\ref{eq:1}) and (\ref{eq:2}) are used as models in a broad range of applications in science and engineering, in fluid flows in porous media, fracture dynamics, thin film dynamics, surface growth dynamics (e.g. surface erosion by ion sputtering processes), and so on \cite{Alava2004imbibition,Cuerno1995stochastic,Hu2009stochastic,Soriano2005anomalous}, and references cited therein. Thus, it is important to obtain their solutions. There are few works focusing on obtaining solutions to equation (\ref{eq:2}) , see for example \cite{Gao2018numerical}. However, to my knowledge, there is no work that focuses on solving (\ref{eq:1}).  

Deriving analytical solutions for the nonlinear SPDEs  (\ref{eq:1}) and (\ref{eq:2}) is not an easy task at all. Therefore, numerical methods have to be used. In this paper, a numerical method based on Wiener chaos expansion (WCE) is used. Most importantly, WCE based numerical methods separate deterministic effects from random effects in a effective manner, and allow to obtain a sytem of deterministic partial differential equations (PDEs) for the coefficients of the WCE, see for example \cite{Mikulevicius2004stochastic}. Methods based on the WCE have been used to numerically solve stochastic evolution equations such as stochastic Burgers' equation and stochastic Navier-Stokes equation, stochastic vorticity equation \cite{Hou2006Wiener,Nijimbere2014ionospheric}, to name few.  

In the numerical simulations of SPDEs, analytical solutions are importantly needed in the assessment of the accuracy of the numerical solutions. For the linearized SgKS equation
\begin{equation}
 \hspace{.12cm}\frac{\partial u}{\partial t}=-\kappa \frac{\partial^2 u}{\partial x^2}-\eta \frac{\partial^3 u}{\partial x^3}-\nu\frac{\partial^4 u}{\partial x^4}+\sigma \dot{W},
\label{eq:3}
\end{equation}
for example, analytical solutions can be obtained and shall be used in our preliminary tests in section \ref{sec:3}. In the simulations of the quasi-nonlinear SPDEs (\ref{eq:1}) and (\ref{eq:2}), the accuracy of the numerical solutions is evaluated by comparing the numerical solutions to semi-analytical solutions obtained by performing some change of variables that transform the SPDE into a deterministic equation for which a numerical solution can fairly be obtained using an appropriate numerical method.

WCE approximate numerical solutions of the SPDES (\ref{eq:1}) and (\ref{eq:2}) are obtained using  60 terms in the WCE over a time interval up to  $t=3$. The relative difference of the results from the corresponding semi-analytical solution remains small than the desired tolerance level. The  WCE numerical computations took 942.5 sec CPU time compare with 24.0 sec CPU time for the semi-analytical solution on a Pentium (R) PC with 2.5 GHz CPU. However,  the semi-analytical solution required a considerably larger allocation of computer memory. 

The paper is organized as following. A brief description of the WCE method is given in section \ref{sec:2}. In section \ref{sec:3}, the WCE-based method is applied to numerically solve some initial boundary value problems (IBVPs) involving the linearized SPDE in (\ref{eq:3}), and analytical solutions are obtained and compared with the numerical solutions. In section \ref{sec:4}, a semi-analytical solution procedure for the nonlinear SPDEs  (\ref{eq:1}) and  (\ref{eq:2}) is described, and the WCE method is applied to some IBVPs involving these nonlinear  SPDEs. In section  \ref{sec:5}, the results of numerical computations of the IBVPs involving the nonlinear SPDEs (\ref{eq:1}) and (\ref{eq:2}) are presented, and numerical solutions are contrasted with semi-analytical solutions. 

\section{Wiener Chaos Expansion (WCE) method} \label{sec:2}
\begin{definition}
The Wiener chaos expansion (WCE) of a function $u(x,t;{W}_0^t)$ is the infinite series 
\begin{equation}
u(x,t;{W}_0^t)=\sum_\alpha u_\alpha(x,t)T_\alpha,   
\label{eq:4}
\end{equation}
where $u_\alpha(x,t)$ are deterministic functions and $u_\alpha(x,t)=E[ u_\alpha(x,t)T_\alpha]$ is the mean with respect to the noise $W$, and $T_\alpha$ are multivariable Hermite polynomials of Gaussian random variables \cite{Cameron1947orthogonal}.
\label{df:1}
\end{definition}
In the present paper, we consider some IBVPs for which (\ref{eq:4}) is a solution of the SPDE
\begin{equation}
\frac{\partial u}{\partial t}(x,t)=\mathcal{L}[u(x,t)]+\sigma(x)\dot{W}(t),
\label{eq:5}
\end{equation}
where $\mathcal{L}$ is an elliptic differential operator (linear or nonlinear), $\sigma$ is a scalar function of position $x$. In the case of the SgKS equation,
\begin{equation}
\mathcal{L}[u(x,t)]=-u \frac{\partial u}{\partial x}(x,t)-\kappa \frac{\partial^2 u}{\partial x^2}(x,t)-\eta \frac{\partial^3 u}{\partial x^3}(x,t)-\nu\frac{\partial^4 u}{\partial x^4}(x,t).
\label{eq:6}
\end{equation}
Let us now define a set of orthonormal bases $\left\{m_i(s)\right\}_{i=1}^{n},\hspace{.12cm} s>0 \hspace{.12cm} \text{in the Hilbert space} \hspace{.12cm} L^2\{[0,t]\}$, and the random variables $\xi_i$,  $i=1,2,\ldots$, 
\begin{equation}
\xi_i=\int\limits_{0}^{t}m_i(s)\dot{W}(s)ds=\int\limits_{0}^{t}m_i(s)dW(s), \hspace{.25cm} i=1,2,\cdots.
\label{eq:7}
\end{equation} 
In that case,  Brownian motion can be written as a L$\acute{e}$vy-Ciesielski series given by 
\begin{equation}
W(s)=\sum\limits_{i=1}^{\infty}\xi_i\int\limits_{0}^{s}m_i(\tau)d\tau, \hspace{.25cm} 0\leq s\leq t.
\label{eq:8}
\end{equation}
This series converges uniformly for $\forall s\leq t$ and in the mean square sense,
\begin{equation}
E\left[W(s)-\sum\limits_{i=1}^{N}\xi_i\int\limits_{0}^{s}m_i(\tau)d\tau\right]^2\rightarrow 0 \hspace{.15cm}\text{as}\hspace{.15cm} N\rightarrow\infty.
\label{eq:9}
\end{equation} 

In particular, defining the orthonormal basis in terms of trigonometric functions as $m_1(t)={1}/{\sqrt{T}}, m_i(t)=\sqrt{{2}/{T}}\cos\left[(i-1)\pi t/{T}\right],\,\, i=2,3,\cdot\cdot\cdot, \,\,0\leq t\leq T$, gives  the Paley-Wiener representation of $W(t)$ (see, for example,  equation (1.23) in \cite{Mikosch2004elementary})
\begin{equation}
W(t)=\frac{t}{\sqrt{T}}+\frac{\sqrt{2T}}{\pi}\sum\limits_{i=2}^\infty \frac{\xi_i}{i-1} {\sin\left[\frac{(i-1)\pi t}{T}\right]}.
\label{eq:10}
\end{equation}

We can now write solutions to (\ref{eq:5}) as
\begin{equation}
u(x,s;{W}(s))=u(x,s;\xi_1,\cdots,\xi_n,\cdots), 0\le s\le t, 
\label{eq:11}
\end{equation}
and hence, $u$ can be expressed as in (\ref{eq:4}), where 
\begin{equation}
T_{\alpha}(\xi)=\prod\limits_{i=1}^{\infty}H_{\alpha_i}(\xi_i),
\label{eq:12}
\end{equation}
the functions $H_{\alpha_i}(\xi_i)$ are Hermite polynomials of order $\alpha_i$ and are normalized with respect to Gaussian measure, and  $T_{\alpha}(\xi)$ are called Wick polynomials, with $\alpha_i$ defined within the set
\begin{equation}
\mathcal{G}=\left\{\alpha=(\alpha_i,i\geq 1)|\alpha_i\in \{0,1,2,3,\cdot\cdot\cdot\},|\alpha|=\sum\limits_{1}^{\infty}\alpha_i<\infty\right\}.
\label{eq:13}
\end{equation}

Thus, an appropriate construction of the set $\mathcal{G}$ (the choice the $\alpha_i$ values) plays an important role in the convergence and the accuracy of the numerical solution.
We also note that the WCE based numerical methods can be generalized to SPDEs with the stochastic forcing of the form $\sum_{j=1}^{jmax}\sigma_j W_j(t), jmax<\infty$ (see for example \cite{Mikulevicius2004stochastic}).

The parabolocity of (\ref{eq:4}) and the regularity of the noise $W(t)$ imply that the solution of (\ref{eq:4}) is square-integrable with respect to the noise. In that case, the mean and the variance can be expressed in terms of $u_\alpha$ as stated in the following theorem. 

\begin{theorem}{[\textbf{Cameron-Martin}]}
Suppose that for any $x\in\mathbb{R}$ and $s\leq t \in \mathbb{R}$, the solution of $u(x,s)$ of equation (\ref{eq:4}) is a functional of the Brownian motion $\{W(t), 0\leq t\leq T\}$ with $E|u(x,s)|^2<\infty$. Then $u(x,s)$ has the following WCE:
\begin{equation}
u(x,s)=\sum_{\alpha\in\mathcal{G}}u_\alpha(x,s)T_\alpha(\xi), \,\,\,\, u_\alpha(x,s)=E[u(x,s)T_\alpha(\xi)],
\label{eq:14}
\end{equation}
where $T_\alpha(\xi)$ are Wick polynomials defined by equation (\ref{eq:12}), and the mean and variance of $u(x,s)$ are given respectively by
\begin{equation}
E[u(x,s)]=u_0(x,s) \,\,\,\, \text{and} \,\,\,\, E[(u(x,s)-u_0(x,s))^2]=\sum_{\alpha\in\mathcal{G},\alpha \neq 0}|u_\alpha(x,s)|^2
\label{eq:15}
\end{equation}
\label{thm:1}
\end{theorem}
The proof of this theorem can be found in Cameron and Martin \cite{Cameron1947orthogonal}. 

The WCE of a product of two functions is given by the following theorem,
\begin{theorem}
Suppose the functions $u$, $v$ have Wiener chaos expansions
$$u(x,t)=\sum_{\alpha}u_{\alpha}(x,t)T_{\alpha}(\xi), \hspace{.25cm} v(x,t)=\sum_{\beta}v_{\beta}(x,t)T_{\beta}(\xi).$$
If $E(|uv|^2)<\infty$, then the product $uv$ has the Wiener chaos expansion
\begin{equation}
uv=\sum_{\theta\in\mathcal{G}}\left(\sum_{p\in\mathcal{G}}\sum_{0<\beta<\theta}C(\theta,\beta,p)u_{\theta-\beta+p} v_{\beta+p}\right)T_{\theta}(\xi),
\label{eq:16}
\end{equation}
where
$C(\theta,\beta,p)=\left(C_\beta^\theta C_p^{\theta-\beta+p} C_p^{\beta+p}\right)^{1/2}$, with $C_a^b=\frac{b!}{a!(b-a)!}$.
\label{thm:2}
\end{theorem}
This theorem is useful in the evaluation of nonlinear terms, its proof is standard and can be found , for example, in \cite{Mikulevicius2004stochastic} (Lemma 14). 
%
%The expansion for the product $uv$ given in (\ref{uv1}) is analogous to the expression for the convolution in the context of the Fourier transform.

Numerical methods based on the WCE consist of writing the solution of the SPDE in terms of (\ref{eq:4}). This allows to obtain a system of deterministic partial differential equations (PDEs) for the coefficients $u_\alpha$ known as the propagator associated with the SPDE. An appropriate numerical approximation method is then used to numerically solve the propagator. This procedure  is  quite standard. In section \ref{sec:3}, it is applied to obtain the numerical approximate solutions to some IBVPs involving the linearized SgKS equation (\ref{eq:3}), while in section \ref{sec:4}, it is utilized to numerically solve some IBVPs involving the SgKs equation (\ref{eq:1}) and the SKS equation (\ref{eq:2}). 

Numerical implementations are described in  sections \ref{subsec:3.1} and \ref{subsec:4.1}.  A predictor-corrector method, used in \cite{Nijimbere2014ionospheric,Nijimbere2016nonlinear}, which is based on the second order Adam-Bashforth explicit scheme and a third order Adam-Moulton implicit scheme is implemented in order to achieve a fast convergence.

\section{A preliminary test: the WCE method applied to the linearised  stochastic Kuramoto--Sivashinsky equations} \label{sec:3}
In this section, numerical solutions of a test problem involving the linearized SgKS equation (\ref{eq:3}) are obtained.  In the assessment of the accuracy of the numerical solutions, the numerical results are compared with exact solutions of the respective equations. 

Now, consider that the differential operator $\mathcal{L}$ in (\ref{eq:5}) is linear. In that case,  the Wick polynomials are $T_{\alpha_i}=H_{\alpha_i=1}(\xi_i)=\xi_i$ ($i=1, 2, \dots$) since there are no nonlinear product terms in the equation. Thus, the WCE solution of (\ref{eq:5}) becomes  $u=\sum_i u_i\xi_i$, where $\xi_i$ are given by (\ref{eq:7}).
% In that case, integration of (\ref{sto3}) with respect to $t$ yields 

\subsection{Application of the WCE}
\label{subsec:3.1}
Considering that $\mathcal{L}$ is linear and integrating (\ref{eq:5}) with respect to $t$ gives 
\begin{equation}
u(x,t)=u(0,t)+\int\limits_0^t\mathcal{L}[u(x,\tau)]d\tau+\sigma(x)\sum\limits_{i=1}^\infty\xi_i\int\limits_0^t m_i(\tau)d\tau.
\label{eq:17}
\end{equation}
Now, Multiplying both sides of (\ref{eq:17}) by $\xi_i$ and taking the expectation of the resulting equation, while using the fact that the variables $\xi_i$ are independent, we obtain the linear PDE
\begin{equation}
\frac{\partial u_i}{\partial t}(x,t)=\mathcal{L}[u_i](x,t)+\sigma(x) m_i(t),
\label{eq:18}
\end{equation}
which is the propagator associated to the the linear SPDE. For some specified linear operator $\mathcal{L}$,  the propagator  (\ref{eq:18}) can then be solved numerically subject to the appropriate initial and boundary conditions.

Next, let us apply this numerical procedure to some IBVP involving the linearized SgKS equation for which an analytical solution to compare with can be obtained. We consider the IBVP:
\begin{equation}
\frac{\partial u}{\partial t}(x,t)=-\kappa \frac{\partial^2 u}{\partial x^2}(x,t)-\eta \frac{\partial^3 u}{\partial x^3}(x,t)-\nu\frac{\partial^4 u}{\partial x^4}(x,t)+\exp(ikx)\dot{W}(t), \hspace{0.35cm}  t \in (0, \infty), \hspace{0.2cm} x\in(0,2\pi)
\label{eq:19}
\end{equation}
subject to the initial condition 
\begin{equation}
 u(x,0)=V_0 \exp(ikx), \,\,\, V_0 \in \mathbb{R}, \,\,\, k\in\mathbb{N},\,\,\, x\in [0,2\pi], 
\label{eq:20}
\end{equation}
and with periodic boundary conditions.
%\begin{equation}
%u(0,t)=u(2\pi,t), \hspace{0.35cm} t \in [0, \infty).
%\label{SA-D11}
%\end{equation}

We now write its solution as 
\begin{equation}
u(x,t)=V(t)\exp(ikx), 
\label{eq:21}
\end{equation}
substitute into (\ref{eq:20}) and the Langevin  equation
\begin{equation}
d V(t)=(\kappa k^2 +i\eta k^3-\nu k^4)V(t) dt+dW(t), \,\,\, t \in [0, \infty),
\label{eq:22}
\end{equation}
with the initial condition
\begin{equation}
 V(0)=V_0.
\label{eq:23}
\end{equation}
The solution to (\ref{eq:22})-(\ref{eq:23}) is
\begin{equation}
V(t)=\exp\left[(\kappa k^2 +i\eta k^3-\nu k^4)t\right]\left\{V_0+\int\limits_0^t\exp\left[-(\kappa k^2 +i\eta k^3-\nu k^4)\tau\right]dW(\tau)\right\}.
\label{eq:24}
\end{equation}
Hence,
\begin{multline}
u(x,t)=V(t)\exp(ikx)\\=\exp\left[ikx+(\kappa k^2 +i\eta k^3-\nu k^4)t\right]\left\{V_0+\int\limits_0^t\exp\left[-(\kappa k^2 +i\eta k^3-\nu k^4)\tau\right]dW(\tau)\right\}.
\label{eq:25}
\end{multline}

A WCE expression for  (\ref{eq:25}) can be derived. We first observe that the WCE analytical solution of the initial value problem (IVP) (\ref{eq:22})-(\ref{eq:23}) is given by 
\begin{equation}
V(t)=V_0\exp\left[(\kappa k^2 +i\eta k^3-\nu k^4)t\right]+\sum\limits_{i=1}^{\infty}V_i(t)\xi_i,
\label{eq:26}
\end{equation}
where the coefficients of the WCE are $$V_i(t)=\exp\left[(\kappa k^2 +i\eta k^3-\nu k^4)t\right]\left\{\int\limits_0^t\exp\left[-(\kappa k^2 +i\eta k^3-\nu k^4)\tau\right]m_i(s)ds\right\}, i\ge1,$$ 
and do satisfiy the propagator of Langevin equation (\ref{eq:22})  
\begin{equation}
\frac{d V_i (t)}{dt}=(\kappa k^2 +i\eta k^3-\nu k^4) V_i(t)+m_i(t)
\label{eq:27}
\end{equation}
with initial conditions
\begin{equation}
V_0(0)=V_0, \,\,\, V_i(0)=0, i\ge1.
\label{eq:28}
\end{equation}
The WCE analytical solution is therefore
\begin{align}
u_W(x,t)&=\exp(ikx)V_W(t)=\exp(ikx)\sum\limits_{i=0}^\infty V_i\xi_i\nonumber \\&=V_0\exp\left[(ikx+\kappa k^2 +i\eta k^3-\nu k^4)t\right]+\exp(ikx)\sum\limits_{i=1}^{\infty}V_i(t)\xi_i,
\label{eq:29}
\end{align}
and is a WCE of the solution (\ref{eq:25}).

In the numerical computations, on the other hand, the second-order Adams-Bashforth time-discretization scheme is used to numerically solve the propagator, and computations are performed on the domain $[a,b]\times[0,T]=[0,2\pi]\times[0,3]$.  The white noise is generated using the \emph{randn} Matlab function so that $dW=\mbox{randn}(1,N)\sqrt{dt}$, where $dt$ is  the variance of $dW$.
 
Parameters in the linearized SgKS equation (\ref{eq:19}) are respectively set to $\kappa=0.002, \eta=0.002, \nu=0.005$  $k=1$, and $V_0=1$ so that the amplitude of the random forcing  is given by $\sigma(x)=e^{ix}$, while the initial condition is $u(x,0)=e^{ix}$. The time step is set to $\Delta t=0.005$, this corresponds to a time interval of $1000$ time steps. The computation of the solution takes 14.7 sec CPU time on a Pentium (R) PC with 2.5 GHz CPU.

It is important to point out that solution (\ref{eq:24}) involves a nonlocal integral in time which has to be approximated numerically.  The computational expense is diminished by evaluating the integral from $t=0$ to $t=t_n$ as the sum of two integrals, one from $t=0$ to $t=t_{n-1}$ and the other from $t=t_{n-1}$ to $t=t_n$. Each integral is then numerically evaluated using the trapezoidal rule, see for example by \cite{Nijimbere2016nonlinear}. Indeed, the solution (\ref{eq:24}) is a semi-analytical solution.

\subsection{Accuracy assessment of numerical solutions}
\label{subsec:3.2}
To assess the accuracy of the WCE based numerical solutions, we evaluate the absolute difference between the WCE based numerical solution $u_{W}$ and the semi-analytical solution $u$ on the domain $[a,b]\times[0,T]$. Discretizing $x$ as $x=x_k = k \Delta x$  ($k = 0, 1, \dots, K=(b-a)/\Delta x$) , we evaluate the absolute difference over $[a,b]$ at $t=t_n = n \Delta t$ ($n=  0, 1, \dots, N=T/\Delta t$) as
\begin{multline}
\Delta_a u(t_n)=\frac{1}{K}||u_{\tiny W,I}(x_k,t_n) - u(x_k,t_n)||\\= \frac{1}{K}\sum\limits_{k=0}^K|u_{\tiny W,I}(x_k,t_n) - u(x_k,t_n)| = \frac{1}{K}\sum\limits_{k=0}^K|u_{\tiny W,I}^{k}(t_n) - u^{k}(t_n)|, 
\label{eq:30}
\end{multline}
where $u$ is the solution (\ref{eq:25}), $u_{W,I}=\sum_{i\le I} u_i\xi_i$ is the truncated WCE analytical solution and $I$ is the truncation order, and corresponds to the number of terms retained in the Paley-Wiener series of $W(t)$ (the expansion (\ref{eq:10})) in the numerical implementation. 

It can readily be shown that the order of convergence of the WCE based numerical computations can be approximated as 
\begin{equation}
|u_{\tiny W,I}^{k,n} - u^{k,n}|=\sqrt{\Delta t}\sum \limits_{i>I}||V_i||={ \mathcal C}(T,I) (\Delta t)^{p+\frac{1}{2}},
\label{eq:31}
\end{equation}
where $ (\Delta t)^{p}$ is the order of convergence in the numerical implementation of $V_i$, ${\mathcal C}$ is constant that depends on the truncation order $I$ and the length of the time interval $T$, and will become large if $I$ the order of the WCE is small and as the length of the time interval $T$ becomes large. Therefore, the error shall be minimized if a small time step size $\Delta t$ is used in the numerical computations.  

We also evaluate the relative difference between the WCE based numerical solution $u_{W}$ and the analytical solution $u$ over the interval $x \in [a,b]$ for each $t=t_n = n \Delta t$ as
\begin{equation}
\Delta_r u(t_n) =\frac{||u_{\tiny W,I}(x_k,t_n)-u(x_k,t_n)||}{||u_{\tiny W,I}(x_k,t_n)||} = \frac{\sum \limits_{k=0}^{K}|u_{\tiny W,I}^{k,n}(t_n) - u^{k,n}(t_n)|}{\sum \limits_{k=0}^{K}|u^{k,n}(t_n)|}.
\label{eq:32}
\end{equation}

\begin{figure}[!htb]
\begin{picture}(0,0)
\put (20,150){(a)}
\put (220,150){(b)}
\put (20,-20){(c)}
\put (220,-20){(d)}
\put (-25,75){$u(x,t)$}
\put (100,-5){$x$}
\put (235,75){$u(x,t)$}
\put (295,-5){$x$}
\put (-25,-95){$u(x,t)$}
\put (100,-175){$x$}
\put (295,-175){$x$}
\end{picture}
\vspace{.75cm}
\centering{\includegraphics[scale=0.5]{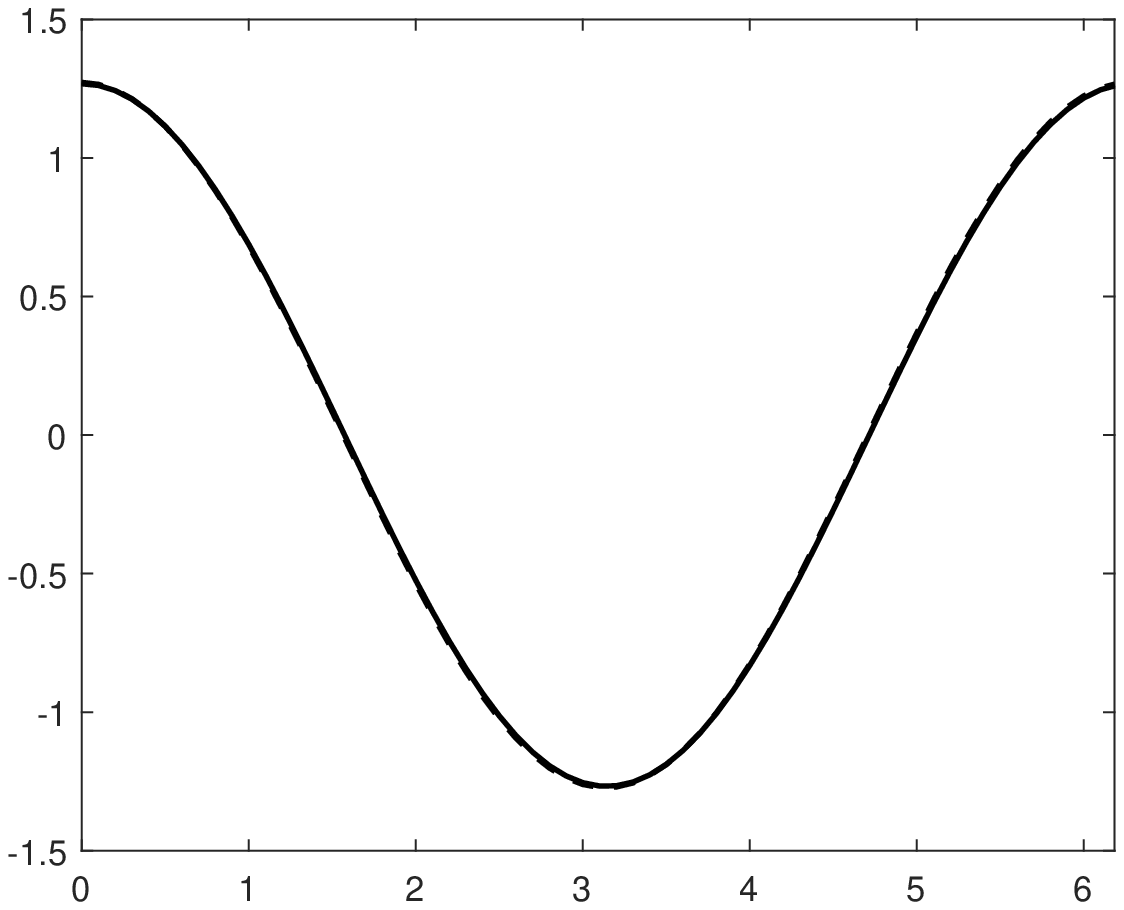}
%\hspace{.75cm}
 \includegraphics[scale=0.5]{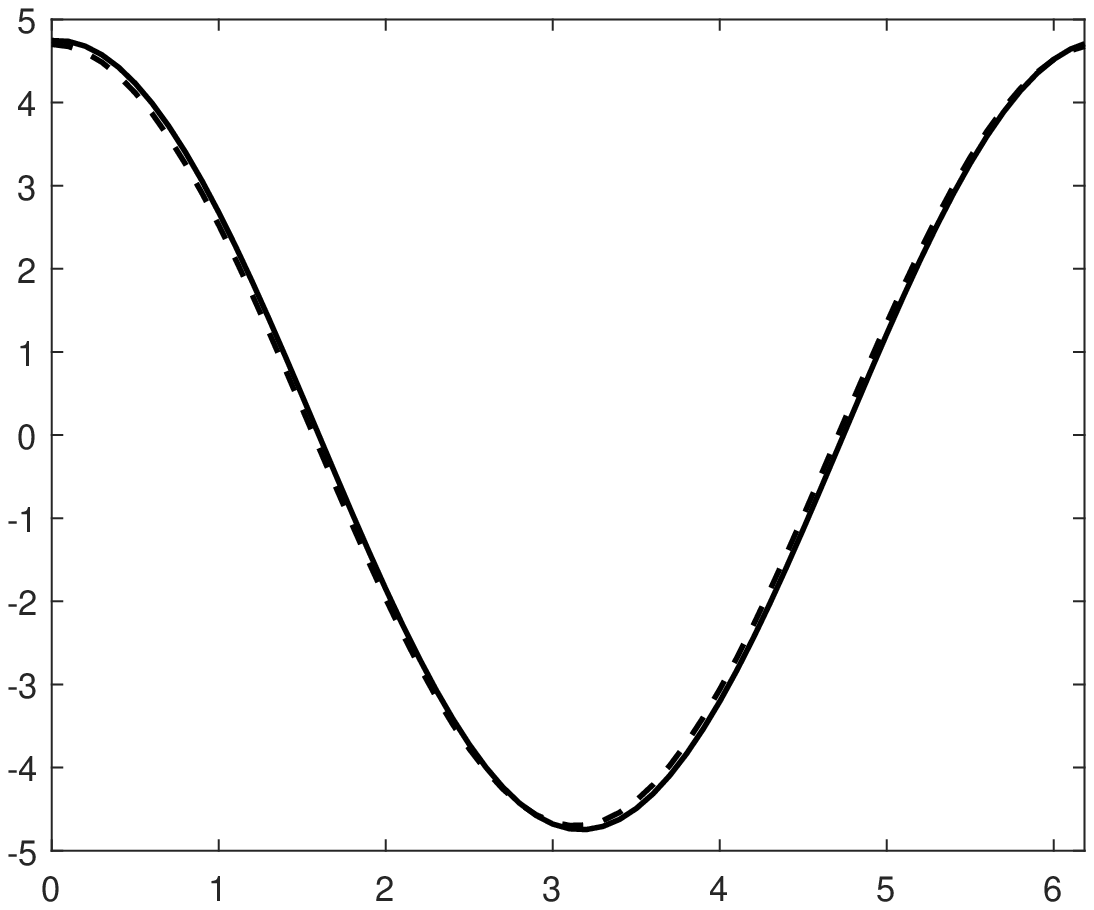}}
\centering{\includegraphics[scale=0.5]{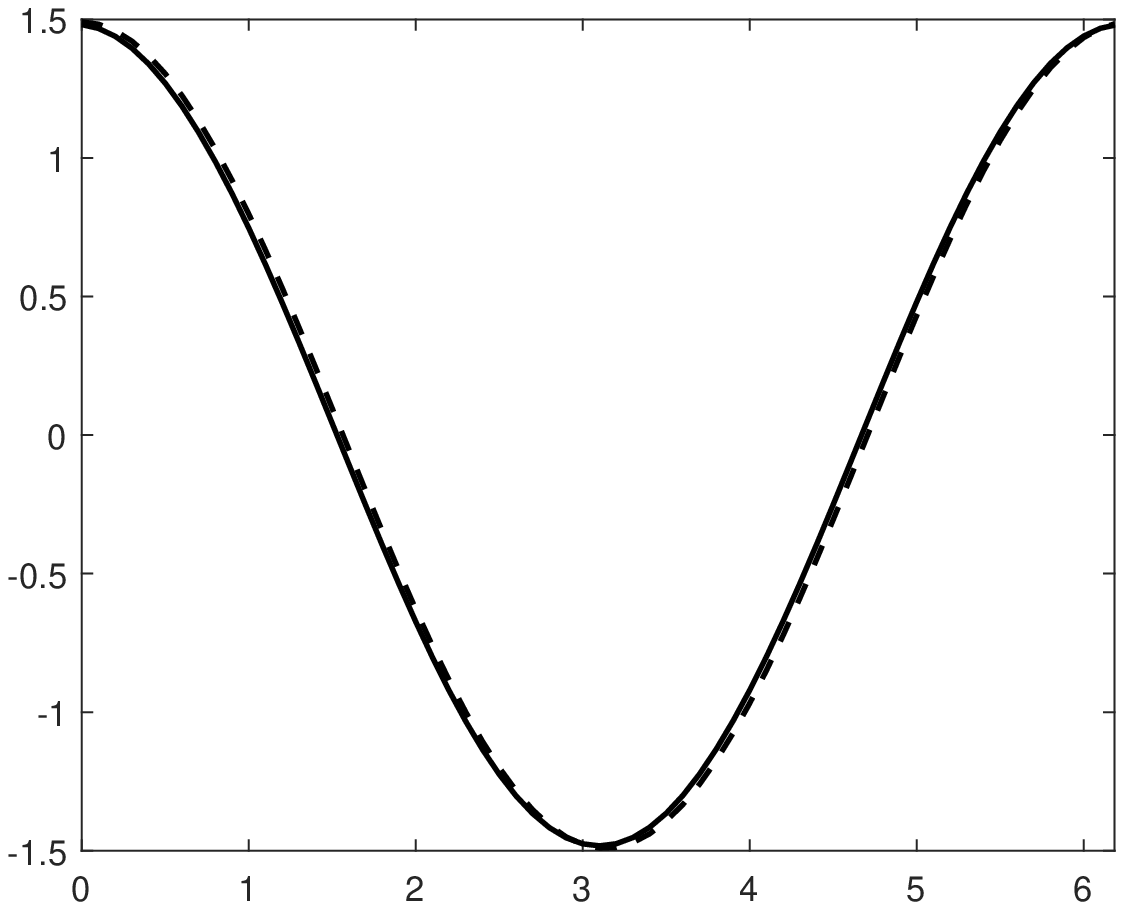}
%\hspace{.75cm} 
\includegraphics[scale=0.5]{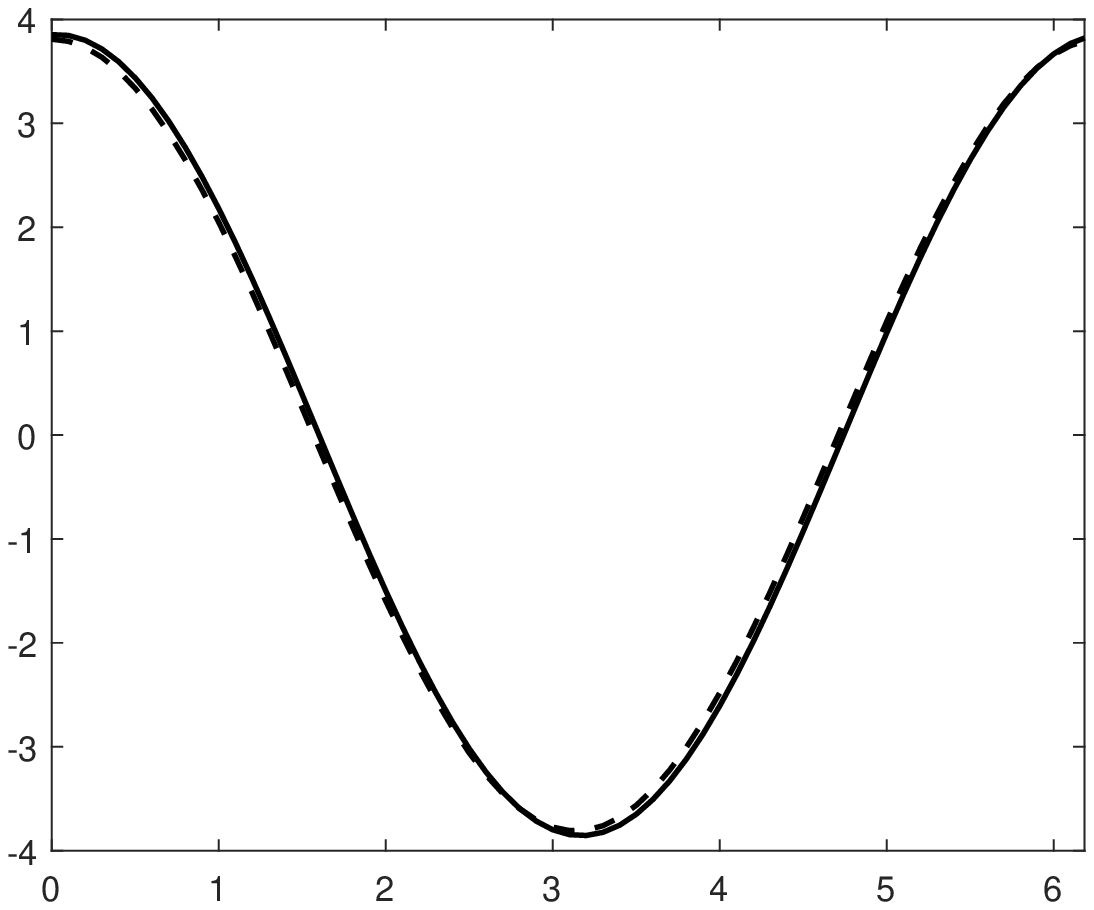}}
\caption{The linearized SgKS equation (\ref{eq:19}): solution $u(x,t)$ as a function of the position $x$ at the time $t=3$. (a) first realization, (b) second realization, (c) third realization and (d) fourth realization. The dashed curve represents the analytical solution (\ref{eq:25}) while the dashed curve represents the WCE based numerical solution.}
\label{fig:1}
\end{figure}

\begin{figure}[!htb]
\begin{picture}(0,0)
\put (20,150){(a)}
\put (220,150){(b)}
\put (20,-20){(c)}
\put (220,-20){(d)}
\put (-25,75){$u(x,t)$}
\put (100,-5){$t$}
\put (295,-5){$t$}
\put (-25,-95){$u(x,t)$}
\put (100,-175){$t$}
\put (295,-175){$t$}
\end{picture}
\vspace{.75cm}
\centering{\includegraphics[scale=0.5]{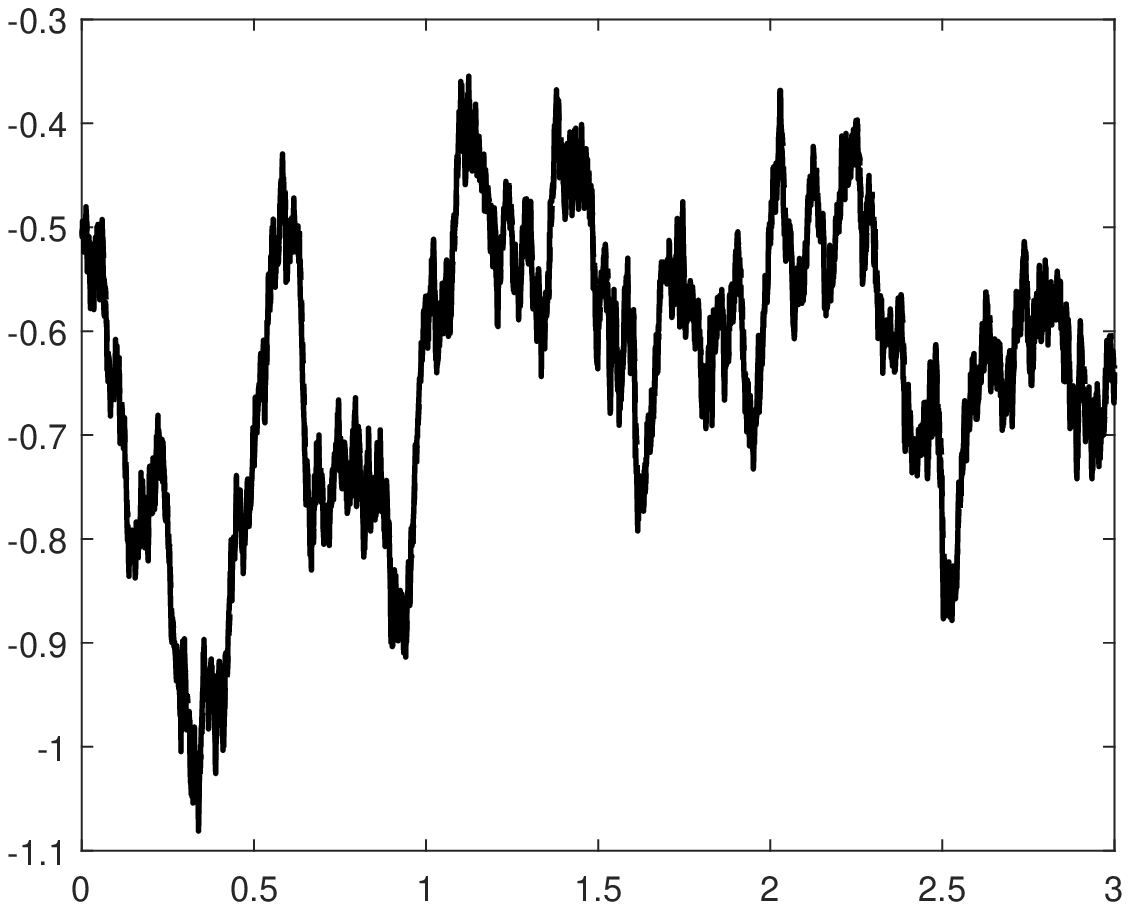} \includegraphics[scale=0.5]{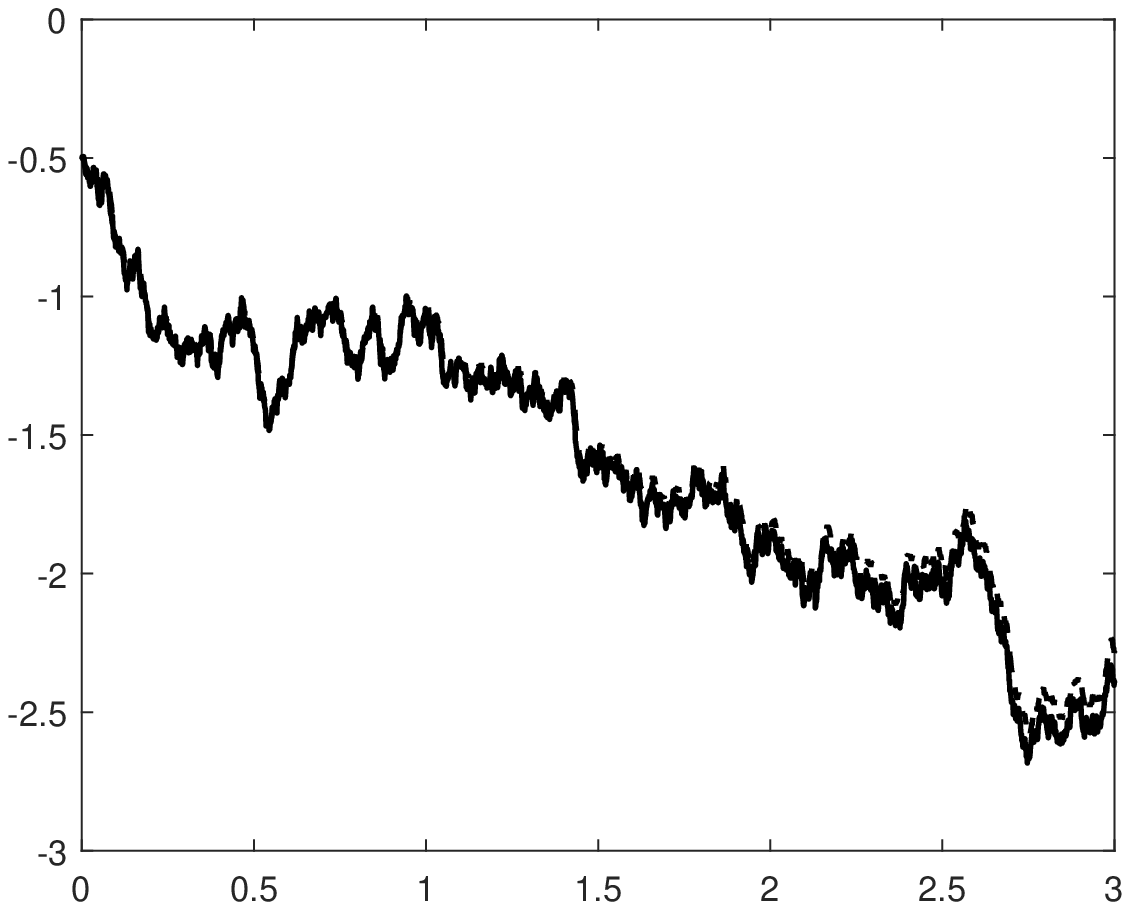}}
\centering{\includegraphics[scale=0.5]{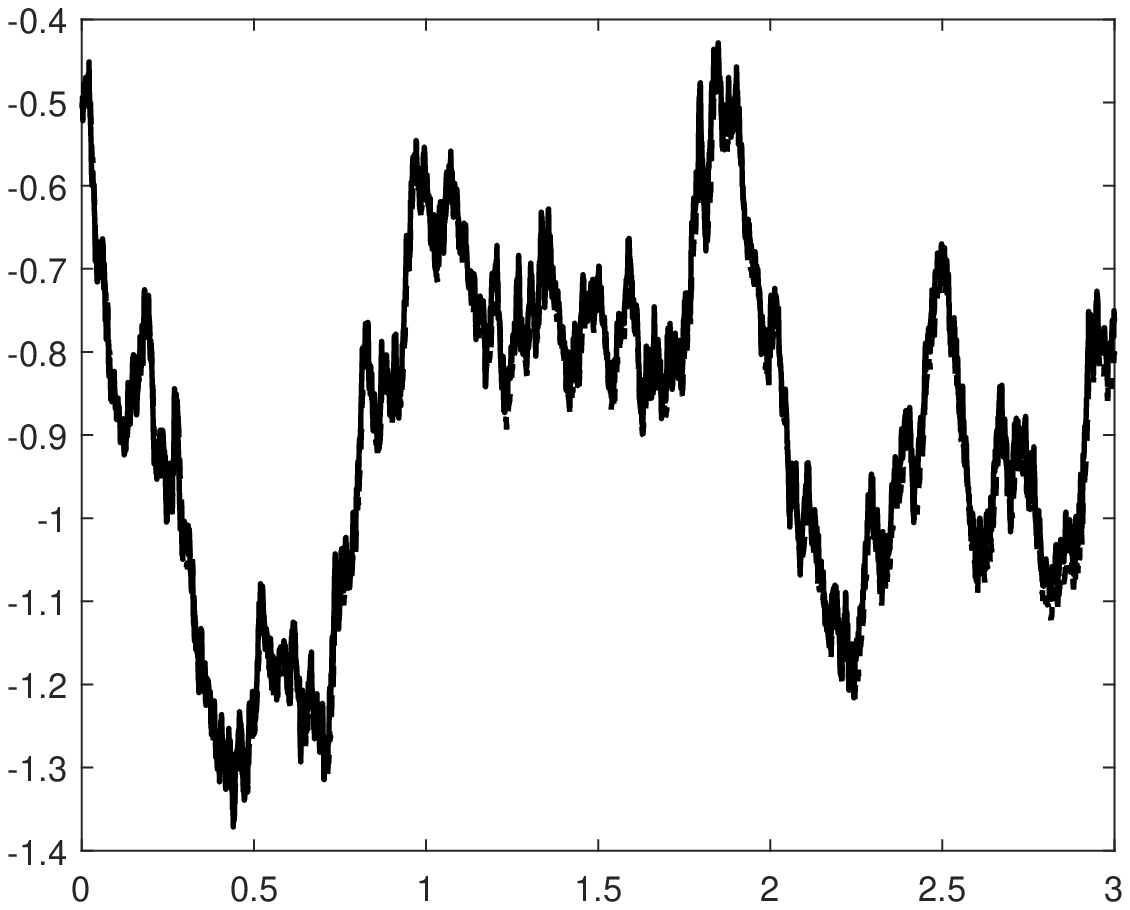} \includegraphics[scale=0.5]{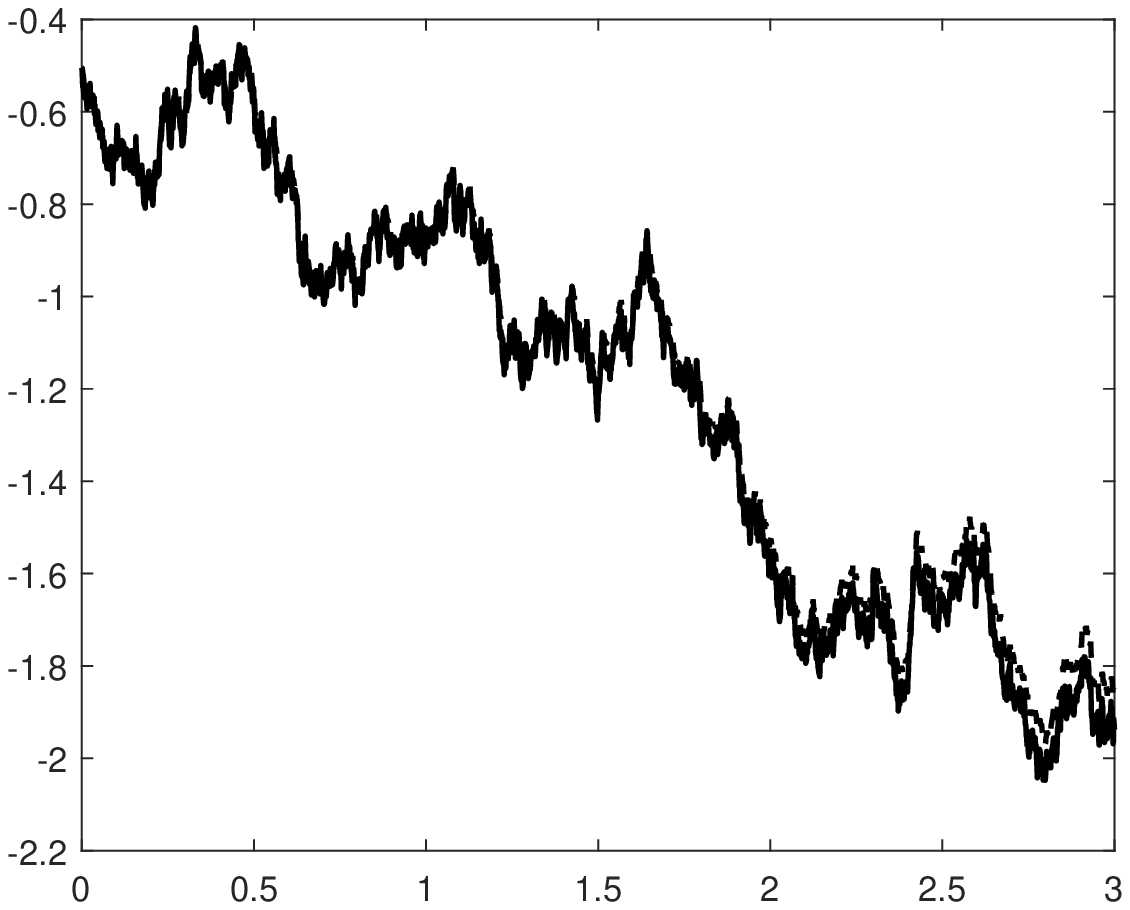}}
\caption{The linearized SgKS equation (\ref{eq:19}): solution $u(x,t)$ as a function of time $t$ at the position $x=1.5$. (a) first realization, (b) second realization, (c) third realization and (d) fourth realization. The dashed curve represents the analytical solution (\ref{eq:25}) while the dashed curve represents the WCE based numerical solution.}
\label{fig:2}
\end{figure}

\begin{figure}[!htb]
\begin{picture}(0,0)
\put (20,150){(a)}
\put (220,150){(b)}
\put (20,-20){(c)}
\put (220,-20){(d)}
\put (-30,75){$\Delta_a u(t)$}
\put (100,-5){$t$}
\put (295,-5){$t$}
\put (-30,-95){$\Delta_a u(t)$}
\put (100,-175){$t$}
\put (295,-175){$t$}
\end{picture}
\vspace{.75cm}
\centering{\includegraphics[scale=0.5]{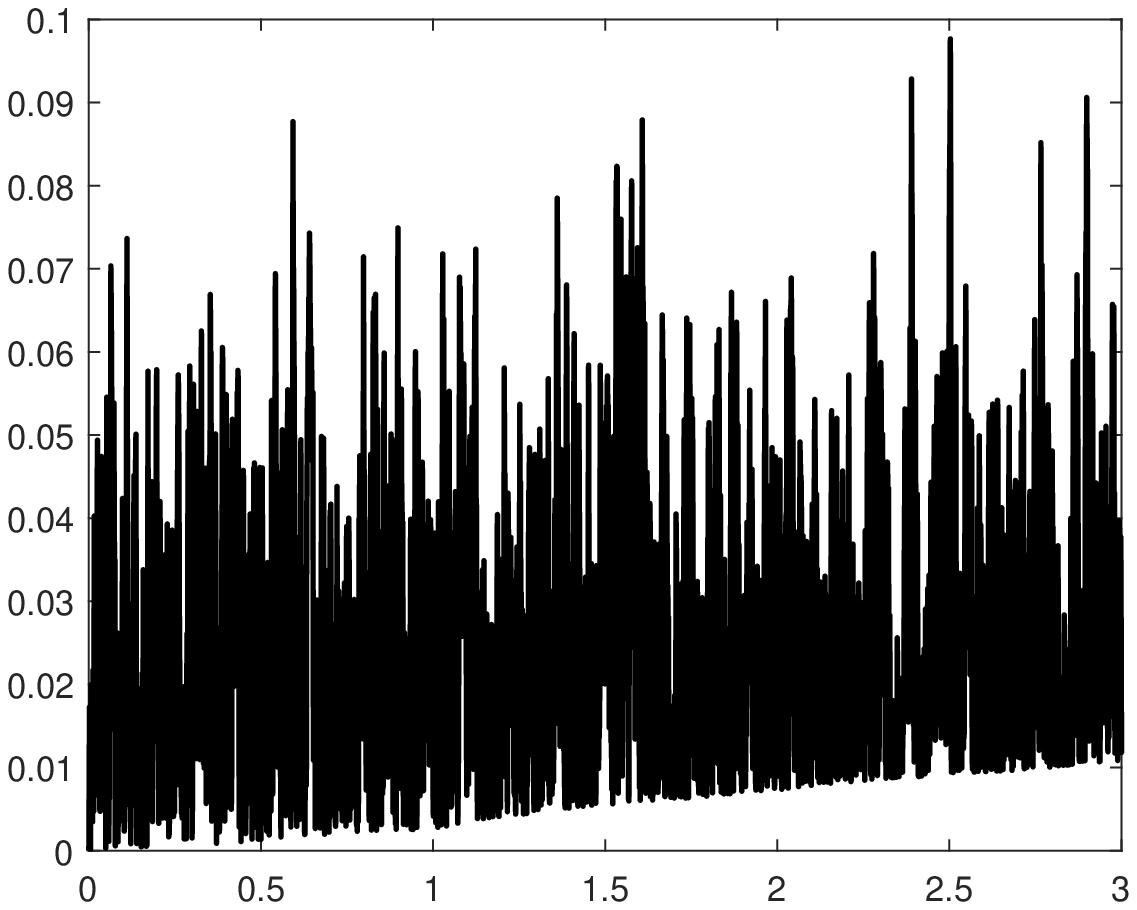} \includegraphics[scale=0.5]{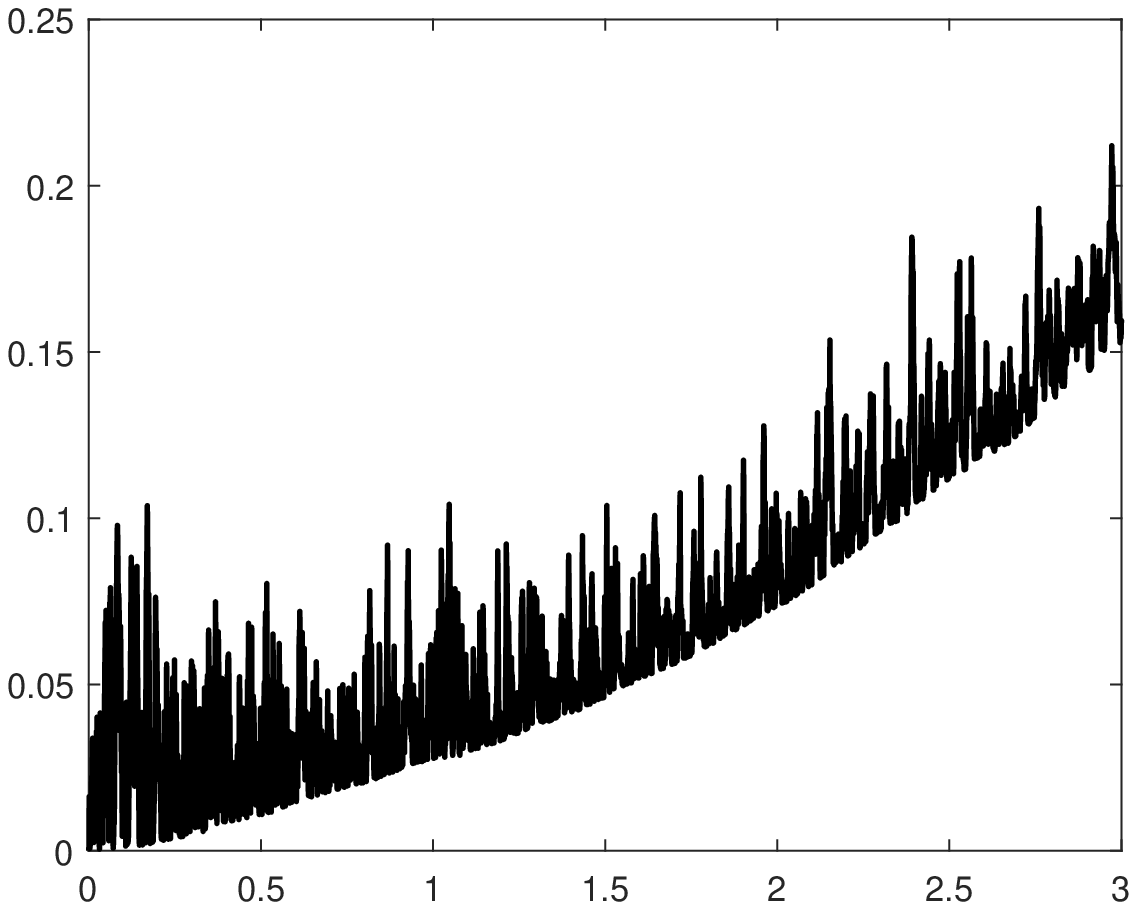}}
\centering{\includegraphics[scale=0.5]{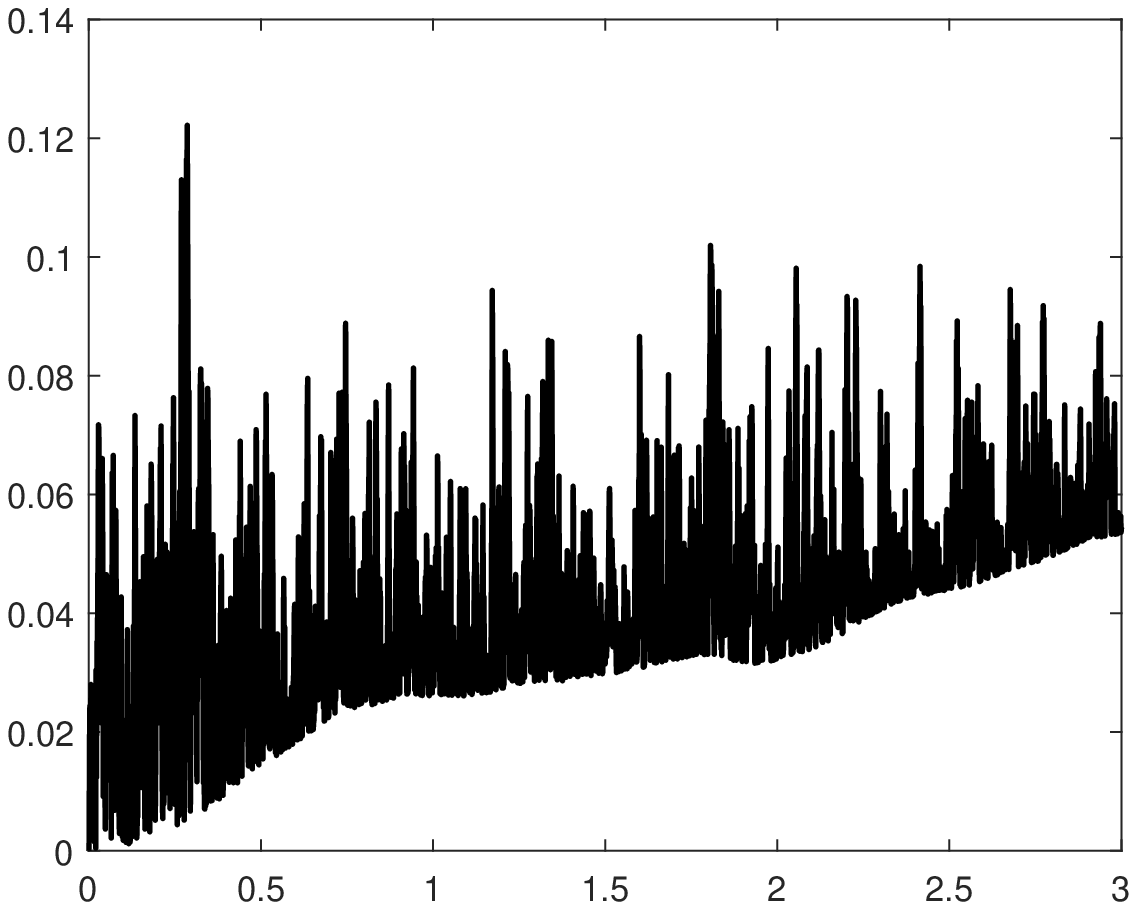} \includegraphics[scale=0.5]{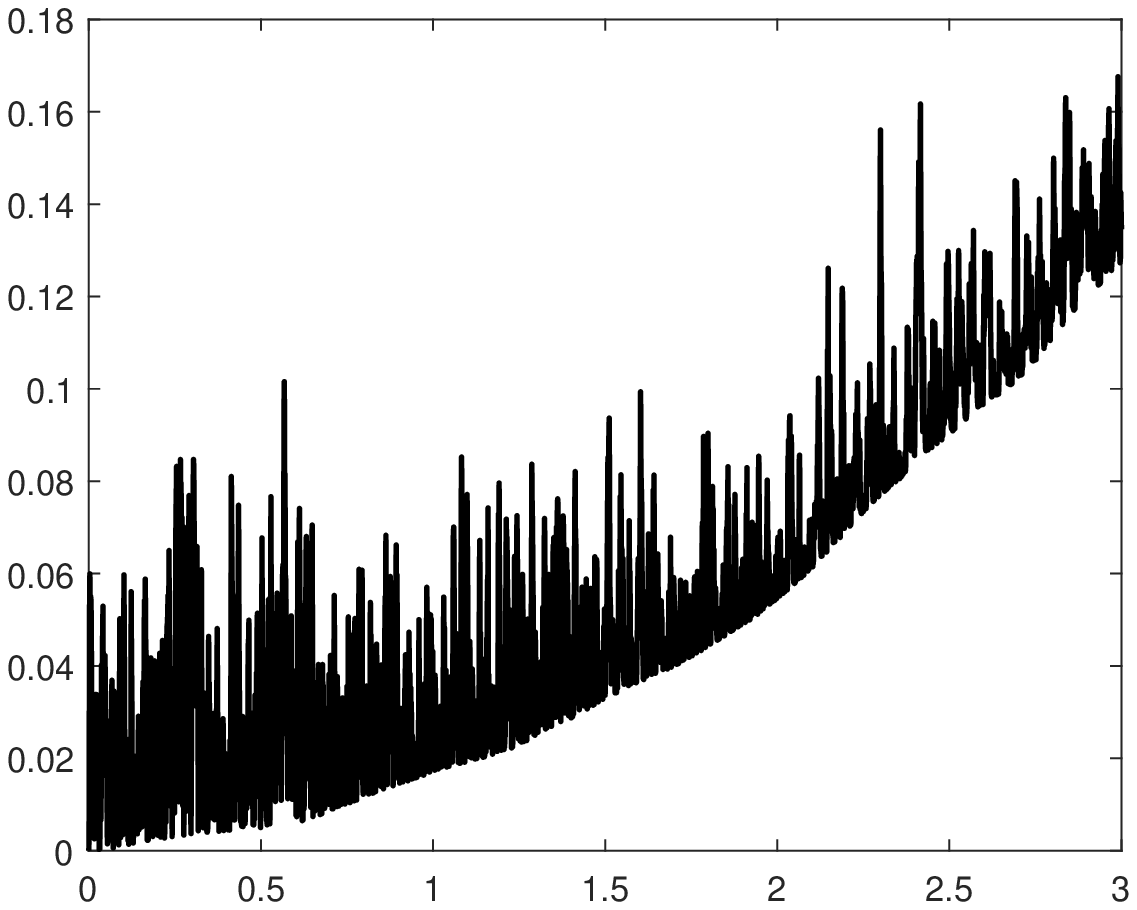}}
\caption{The linearized SgKS equation (\ref{eq:19}): absolute error $\Delta_a u(t)$ (\ref{eq:30}) as a function of time $t$. (a) first realization, (b) second realization, (c) third realization and (d) fourth realization. The dashed curve represents the analytical solution (\ref{eq:25}) while the dashed curve represents the WCE based numerical solution.}
\label{fig:3}
\end{figure}

\begin{figure}[!htb]
\begin{picture}(0,0)
\put (20,150){(a)}
\put (220,150){(b)}
\put (20,-20){(c)}
\put (220,-20){(d)}
\put (-30,75){$\Delta_r u(t)$}
\put (100,-5){$t$}
\put (295,-5){$t$}
\put (-30,-95){$\Delta_r u(t)$}
\put (100,-175){$t$}
\put (295,-175){$t$}
\end{picture}
\vspace{.75cm}
\centering{\includegraphics[scale=0.5]{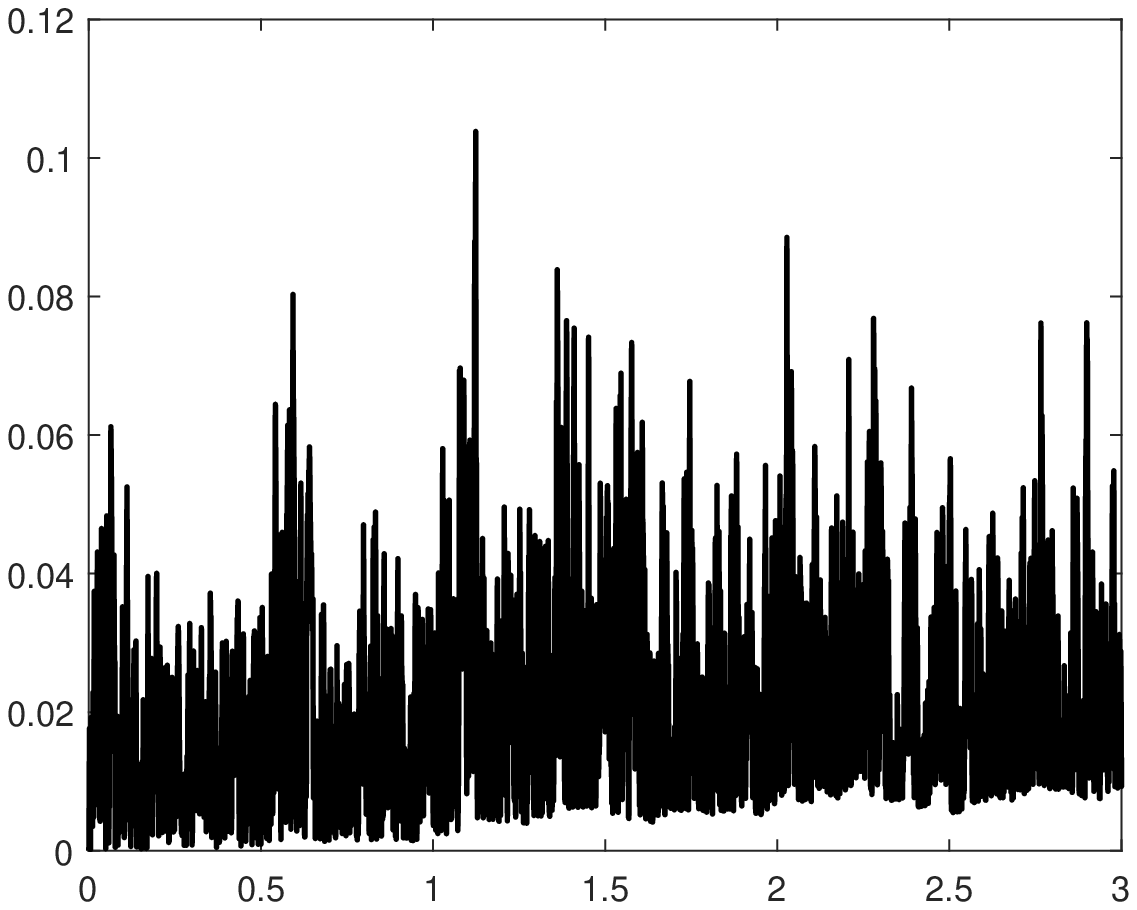} \includegraphics[scale=0.5]{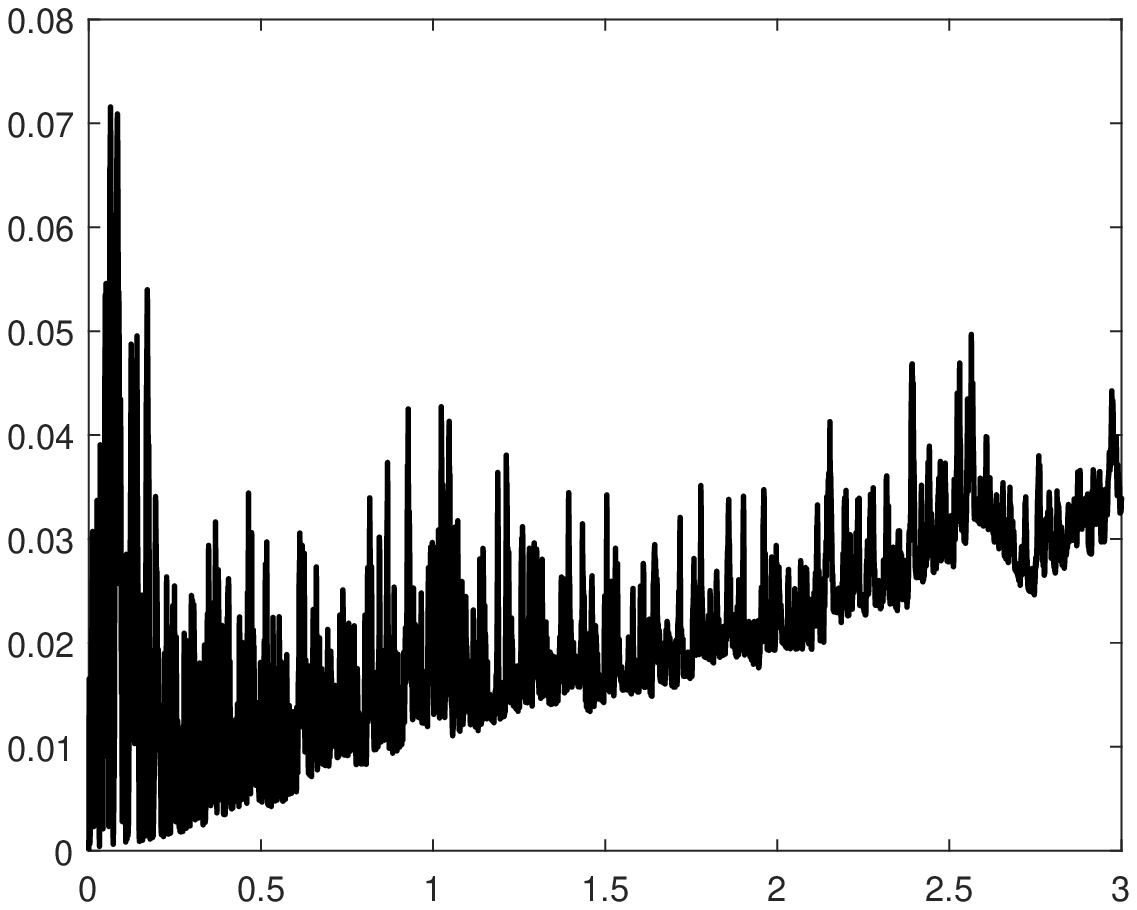}}
\centering{\includegraphics[scale=0.5]{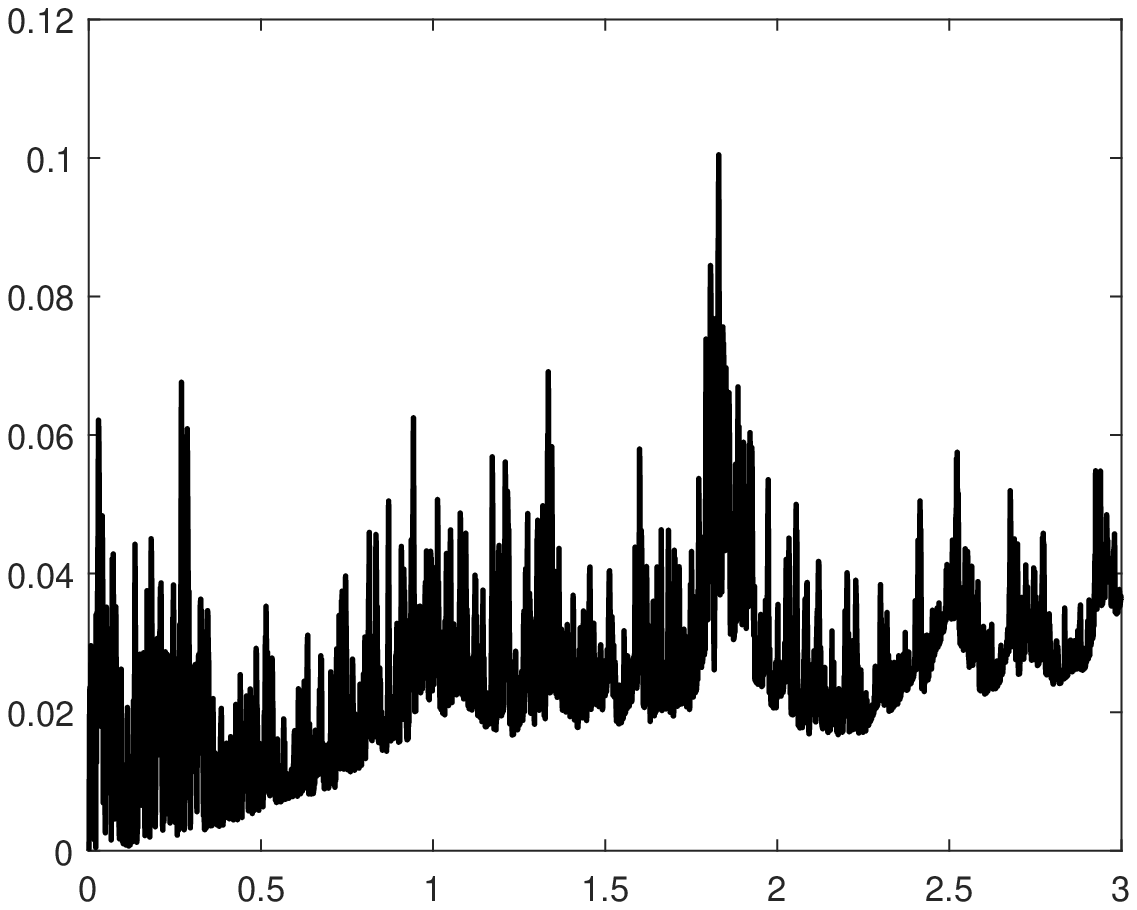} \includegraphics[scale=0.5]{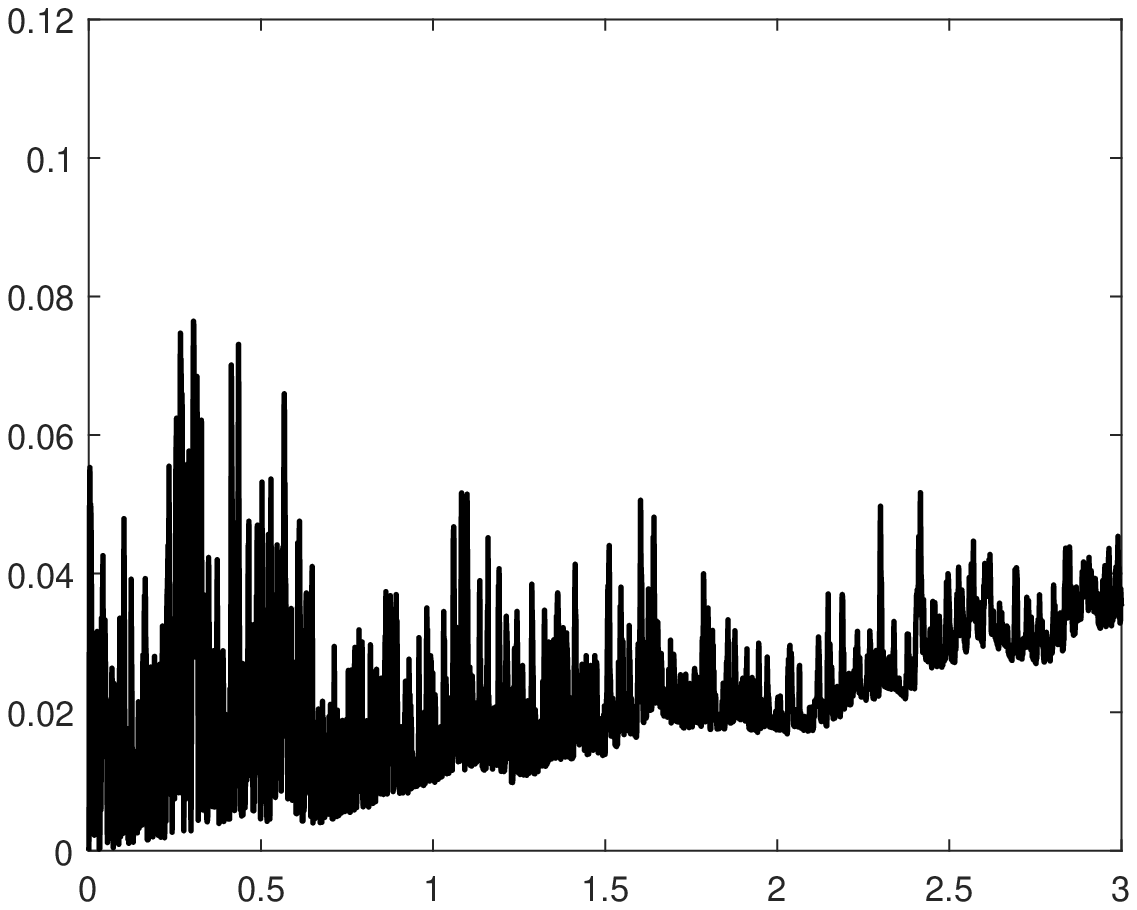}}
\caption{The linearized SgKS equation (\ref{eq:19}): relative error $\Delta_r u(t)$ (\ref{eq:32}) as a function of time $t$. (a) first realization, (b) second realization, (c) third realization and (d) fourth realization. The dashed curve represents the analytical solution (\ref{eq:25}) while the dashed curve represents the WCE based numerical solution.}
\label{fig:4}
\end{figure}

The analytical solutions and the WCE based numerical solutions of the IBVP (\ref{eq:19})-(\ref{eq:20}) involving the linearized SgKS equation are contrasted in Figures \ref{fig:1} to \ref{fig:4}.  In total, four realizations were performed. In each realization, the absolute and relative differences are evaluated.  The WCE based numerical solutions of the linearized SgKS equation as functions of $x$ at the time of $t=3$ are shown in Figure \ref{fig:1}, while they are shown in Figure \ref{fig:2} as functions of time $t$ at $x = 1.5$.  In all four realizations,  there is good agreement between the numerical and analytical solutions over the time interval of $t \in [0,3]$. 

The difference between the analytical and numerical solutions is so small that the dashed curves corresponding to the analytical solution are not clearly visible in Figures \ref{fig:1} and \ref{fig:2}. The absolute difference (\ref{eq:30}) and the relative difference  (\ref{eq:32}) are shown in Figures \ref{fig:3} and \ref{fig:4} respectively. As seen in Figure \ref{fig:3}, the absolute difference (error) increases with time as predicted by equation (\ref{eq:31}).  It is also seen in Figure \ref{fig:4} that the relative difference also increases with time, and is order $10^{-1}$ or less over the time interval $[0,3]$.  Thus, th error can be minimized using a small time step size $\Delta t$ and a higher order WCE as predicted by (\ref{eq:31}) .

\section{Application of the WCE Method to the Stochastic  Generalized Kuramoto--Sivashinsky Equation}\label{sec:4}
\subsection{WCE based numerical implementation} \label{subsec:4.1}

In this section, we focus on the quasi-nonlinear SPDEs (\ref{eq:1}) and (\ref{eq:2}) and describe the WCE based numerical implementation for these equation. We first consider the following IBVP:
\begin{equation}
\frac{\partial u}{\partial t}=-u \frac{\partial u}{\partial x}-\kappa \frac{\partial^2 u}{\partial x^2}-\eta \frac{\partial^3 u}{\partial x^3}-\nu\frac{\partial^4 u}{\partial x^4}+\sigma \dot{W}(t),
\label{eq:33}
\end{equation}
with the initial condition
\begin{equation}
u(x,0)=f(x), \hspace{.35cm} x\in [a,b]
\label{eq:34}
\end{equation}
and the stochastic mixed (Robin) boundary conditions 
\begin{equation}
a_1\hspace{0.05cm} u(a,t)+a_2 \hspace{0.05cm} u_x(a,t)=g_1(t),\hspace{0.2cm} t \in (0, \infty),
\label{eq:35}
\end{equation} 
\begin{equation}
b_1\hspace{0.05cm} u(b,t)+b_2 \hspace{0.05cm}u_x(b,t)=g_2(t), \hspace{0.2cm} t \in (0, \infty),
\label{eq:36}
\end{equation}
\begin{equation}
c_1\hspace{0.05cm} u(a,t)+c_2 \hspace{0.05cm} u_x(a,t)=g_3(t),\hspace{0.2cm} t \in (0, \infty),
\label{eq:37}
\end{equation} 
and
\begin{equation}
d_1\hspace{0.05cm} u(b,t)+d_2 \hspace{0.05cm}u_x(b,t)=g_4(t), \hspace{0.2cm} t \in (0, \infty),
\label{eq:38}
\end{equation}
where $a_1, a_2, b_1,c_1,c_2, d_1$ and $d_2$ are constants and $g_1(t), g_2(t), g_3(t)$ and $g_4(t)$ may be deterministic functions or stochastic functions of $t$.

 To obtain the propagator associated to (\ref{eq:33}) (a system of equations for the WCE coefficients $u_\alpha$), we first write (\ref{eq:33}) in integral form as
\begin{equation}
u(x,t)=u_0(x)-\int\limits_{0}^{t}\left[u(x,\tau) \frac{\partial u}{\partial x}(x,\tau)+\kappa \frac{\partial^2 u}{\partial x^2}(x,\tau)+\eta \frac{\partial^3 u}{\partial x^3}(x,\tau)+\nu\frac{\partial^4 u}{\partial x^4}(x,\tau)\right]d\tau+\sigma(x)W(t).
\label{eq:39}
\end{equation}
Expressing the solution of (\ref{eq:33}) in terms of the WCE as $u=\sum_\alpha u_\alpha T_\alpha$, multiplying both sides by $T_\alpha$ and taking the expectation yields
\begin{multline}
u_\alpha(x,t)=u_0(x)\mathbb{I}_{\alpha=0}-\int\limits_{0}^{t}\Bigl\{E \Bigl[\Bigl(u\frac{\partial u}{\partial x}\Bigr) T_\alpha\Bigr](x,\tau)+\kappa \frac{\partial^2 u_\alpha}{\partial x^2}(x,\tau)\\+\eta \frac{\partial^3 u_\alpha}{\partial x^3}(x,\tau)+\nu\frac{\partial^4 u_\alpha}{\partial x^4}(x,\tau)d\tau\Bigr\}+\sigma(x)E[W(t)T_\alpha],
\label{eq:40}
\end{multline}
where $\mathbb{I}_{\alpha=0} = 1$ if $\alpha=0$ and zero otherwise.

Applying Theorem \ref{thm:2} gives
\begin{equation}
E \left[\left(u\frac{\partial u}{\partial x}\right) T_\alpha\right]=\sum_{p\in\mathcal{G}}\sum_{0\leq\beta\leq\alpha} C(\alpha,\beta, p)u_{\alpha-\beta+p}\frac{\partial u_{\beta+p}}{\partial x}.
\label{eq:41}
\end{equation}
Using the L$\acute{e}$vy-Ciesielski series representation of $W(t)$ in (\ref{eq:8}) yields
\begin{equation}
E[W(t) T_\alpha]=\sum\limits_{i=1}^{\infty}\int\limits_{0}^{t} m_i(\tau)d\tau  E[\xi_i T_\alpha].
\label{eq:42}
\end{equation}
Observing that $\xi_i=H_1(\xi_i)=T_{\alpha_i=\delta_{i,j}}(\xi_i)$ and that Hermite polynomials are orthogonal with respect to the Gaussian measure, we obtain
$$E[\xi_i T_\alpha]=cov\left[H_1(\xi_i),\prod\limits_{i=1}^{\infty} H_{\alpha_i}(\xi_i)\right]=\mathbb{I}_{\alpha_j=\delta_{i,j}},$$
where $\mathbb{I}_{\alpha_j=\delta_{i,j}}=1$ if $\alpha_i=\delta_{i,j}$ and zero otherwise. This gives
\begin{equation}
E[W(t) T_\alpha]=\sum\limits_{i=1}^{\infty}\mathbb{I}_{\alpha_j=\delta_{i,j}}\int\limits_{0}^{t} m_i(\tau)d\tau.
\label{eq:43}
\end{equation}
Substituting (\ref{eq:41}) and (\ref{eq:43}) into (\ref{eq:40}) and differentiating the resulting equation with respect to $t$, we obtain the equation for the coefficients $u_\alpha$ in the WCE $u=\sum_\alpha u_\alpha T_\alpha$ which is 
%\begin{multline}
%u_\alpha(x,t)=u_0(x)\mathbb{I}_{\alpha=0}+\int\limits_{0}^{t}\nu \frac{\partial^2 u_\alpha}{\partial x^2}(x,\tau)d\tau+\sigma(x)\sum\limits_{i=1}^{\infty}\mathbb{I}_{\alpha_j=\delta_{i,j}}\int\limits_{0}^{t} m_i(\tau)d\tau\\-\sum_{p\in\mathcal{G}}\sum_{0\leq\beta\leq\alpha} C(\alpha,\beta, p)\int\limits_{0}^{t}
%u_{\alpha-\beta+p}(x,\tau)\frac{\partial u_{\beta+p}}{\partial x}(x,\tau)d\tau.\nonumber
%\end{multline}
%Differentiating this with respect to $t$ gives the propagator for the stochastic Burgers' equation (\ref{burger1}),
\begin{multline}
\frac{\partial u_\alpha}{\partial t}(x,t)=-\sum_{p\in\mathcal{G}}\sum_{0\leq\beta\leq\alpha} C(\alpha,\beta, p)u_{\alpha-\beta+p}(x,t)\frac{\partial u_{\beta+p}}{\partial x}(x,t)-\kappa \frac{\partial^2 u_\alpha}{\partial x^2}(x,t)\\-\eta \frac{\partial^3 u_\alpha}{\partial x^3}(x,t)-\nu\frac{\partial^4 u_\alpha}{\partial x^4}(x,t)+\sigma(x)\sum\limits_{i=1}^{\infty}\mathbb{I}_{\alpha_j=\delta_{i,j}} m_i(t),
\label{eq:44}
\end{multline}
and where the indices $\alpha$ are defined within the set in the equation (\ref{eq:13}).  

Writing the initial condition (\ref{eq:34}) and the boundary conditions (\ref{eq:35})-(\ref{eq:38}) in the WCE form, 
implies that the WCE coefficients corresponding to $\alpha\neq0$ have to satisfy the homogeneous initial condition
\begin{equation}
 u_{\alpha\neq0}(x,0)=0, \hspace{.35cm} x\in [a,b],
\label{eq:45}
\end{equation}
and the mixed (Robin) boundary conditions
\begin{equation}
a_1\hspace{0.05cm} u_\alpha(a,t)+a_2 \hspace{0.05cm} \frac{\partial u_\alpha}{\partial x}(a,t)=E[g_1(t)T_\alpha], \hspace{0.35cm} t \in (0, \infty),
\label{eq:46}
\end{equation}
\begin{equation}
b_1\hspace{0.05cm} u_\alpha(b,t)+b_2 \hspace{0.05cm}\frac{\partial u_\alpha}{\partial x}(b,t)=E[g_2(t)T_\alpha], \hspace{0.35cm} t \in (0, \infty),
\label{eq:47}
\end{equation}
\begin{equation}
c_1\hspace{0.05cm} u_\alpha(a,t)+c_2 \hspace{0.05cm} \frac{\partial u_\alpha}{\partial x}(a,t)=E[g_3(t)T_\alpha], \hspace{0.35cm} t \in (0, \infty)
\label{eq:48}
\end{equation}
and
\begin{equation}
d_1\hspace{0.05cm} u_\alpha(b,t)+d_2 \hspace{0.05cm}\frac{\partial u_\alpha}{\partial x}(b,t)=E[g_4(t)T_\alpha], \hspace{0.35cm} t \in (0, \infty).
\label{eq:49}
\end{equation}
The WCE coefficient corresponding to $\alpha=0$, on the other hand,  satisfies the deterministic PDE
\begin{equation}
\frac{\partial u_0}{\partial t}(x,t)=-u_0\frac{\partial u_0}{\partial x}-\sum_{\alpha\neq0} u_\alpha\frac{\partial u_\alpha}{\partial x}-\kappa \frac{\partial^2 u_0}{\partial x^2}-\eta \frac{\partial^3 u_0}{\partial x^3}-\nu \frac{\partial^4 u_0}{\partial x^4}
\label{eq:50}
\end{equation}
subject to the deterministic initial condition
\begin{equation}
 u_{0}(x,0)=f(x), \hspace{.35cm} x\in [a,b],
\label{eq:51}
\end{equation}
and the mixed (Robin) boundary conditions
\begin{equation}
a_1\hspace{0.05cm} u_0(a,t)+a_2 \hspace{0.05cm} \frac{\partial u_0}{\partial x}(a,t)=E[g_1(t)], \hspace{0.35cm} t \in (0, \infty),
\label{eq:52}
\end{equation}
\begin{equation}
b_1\hspace{0.05cm} u_0(b,t)+b_2 \hspace{0.05cm}\frac{\partial u_0}{\partial x}(b,t)=E[g_2(t)], \hspace{0.35cm} t \in (0, \infty),
\label{eq:53}
\end{equation}
\begin{equation}
c_1\hspace{0.05cm} u_0(a,t)+c_2 \hspace{0.05cm} \frac{\partial u_0}{\partial x}(a,t)=E[g_3(t)], \hspace{0.35cm} t \in (0, \infty)
\label{eq:54}
\end{equation}
and
\begin{equation}
d_1\hspace{0.05cm} u_0(b,t)+d_2 \hspace{0.05cm}\frac{\partial u_0}{\partial x}(b,t)=E[g_4(t)], \hspace{0.35cm} t \in (0, \infty);
\label{eq:55}
\end{equation}
Thus, we have a system of deterministic equations for the coefficients of the WCE or the propagator.  It is deterministic and can be solved using classical numerical methods. 

In the numerical computations of the WCE coefficients $u_\alpha(x,t)$, a predictor-corrector method is used to numerically solve the propagator.  The dependent variable at a point $(x,t) = (x_k,t_n)  = (a +k \Delta x, n \Delta t), k=0,1\cdots$ and $n=0,1\cdots$, is approximated as
\begin{equation}
{u}_\alpha(x,t) = {u}_\alpha(x_k,t_n)  \approx{u}_\alpha^{k,n}.
\label{eq:56}
\end{equation}
The derivatives of $u_\alpha$ with respect to the independent variable $x$ are numerically evaluated using the second-order central finite difference approximation.

The predictor scheme is made of two time steps. The first time step, uses the Euler method to compute 
$${u}_\alpha^{k,n+1}={u}_\alpha^{k,n}+\Delta t f_\alpha^{k,n},$$
and the explicit second order Adams-Bashforth method is used at subsequent time levels
$${u}_\alpha^{k,n+2}={u}_\alpha^{k,n+1}+\Delta t\left(\frac{3}{2} h_\alpha^{k,n+1}-\frac{1}{2}h_\alpha^{k,n}\right).$$
To improve the accuracy of the solution, the predictor steps are followed by the implicit third-order Adams-Moulton method for the corrector step
$${u}_\alpha^{k,n+2}={u}_\alpha^{k,n+1}+\Delta t\left(\frac{5}{12} h_\alpha^{k,n+2}+\frac{2}{3}h_\alpha^{k,n+1}-\frac{1}{12}h_\alpha^{k,n}\right),$$
where
\begin{multline}
h(x,t,{u}_\alpha(x,t))=-\sum_{p\in\mathcal{G}}\sum_{0\leq\beta\leq\alpha} C(\alpha,\beta, p)u_{\alpha-\beta+p}(x,t)\frac{\partial u_{\beta+p}}{\partial x}(x,t)-\kappa \frac{\partial^2 u_\alpha}{\partial x^2}(x,t)\\-\eta \frac{\partial^3 u_\alpha}{\partial x^3}(x,t)-\nu\frac{\partial^4 u_\alpha}{\partial x^4}(x,t)+\sigma(x)\sum\limits_{i=1}^{\infty}\mathbb{I}_{\alpha_j=\delta_{i,j}} m_i(t).
\label{eq:57}
\end{multline}
In each case, $h(x,t,{u}_\alpha(x,t))$ is approximated by
\begin{multline}
\tilde{h}(x_k, t_n,{u}_\alpha(x_k,t_n), {U}_\alpha(x_k,t_n)) \approx \tilde{h}(x_k, t_n,{u}_\alpha^{k,n}, {U}_\alpha^{k,n}) \approx \tilde{h}_\alpha^{k,n}\\
=-\sum_{p\in\mathcal{G}}\sum_{0\leq\beta\leq\alpha} C(\alpha,\beta, p)u_{\alpha-\beta+p}^{k,n}\frac{u_{\beta+p}^{k,n}-u_{\beta+p}^{k,n}}{2\Delta x}-\kappa {U}_\alpha^{k,n}-\eta \frac{{U}_\alpha^{k+1,n}-{U}_\alpha^{k-1,n}}{2\Delta x}\\ -\nu\frac{{U}_\alpha^{k+1,n}-2{U}_\alpha^{k,n}+{U}_\alpha^{k-1,n}}{\Delta x^2}+\sigma_k\sum\limits_{i=1}^{\infty}\mathbb{I}_{\alpha_j=\delta_{i,j}} m_i^n,
\label{eq:58}
\end{multline}
with
$$U_\alpha(x_k,t_n)=\frac{\partial^2 u_\alpha}{\partial x^2}(x_k,t_n)\approx \frac{{u}_\alpha^{k+1,n}-2{u}_\alpha^{k,n}+{u}_\alpha^{k-1,n}}{\Delta x^2}.$$

The implicit scheme for the corrector step avoids the stiffness of the matrix of coefficients of the finite difference equations. The numerical implementation of this method is simple; and this method allows us to use relatively large time and space increments.

\subsection{A semi-analytical solution procedure for the SgKS equation}\label{subsec:4.2}

Once the WCE based numerical solution has been computed, the accuracy of the numerical solution has to be assessed. However, to my knowledge none has been able to derive an analytical solution to the SgKS equation (\ref{eq:33}) or to the IBVP such as  (\ref{eq:33})- (\ref{eq:38}). By performing a change of variable, the SgKS equation (\ref{eq:33}) can be transformed into a deterministic PDE (the generalized Kuramoto-Sivashinsky equation) with stochastic initial and boundary conditions. The new IBVP can numerically be solved using an appropriate numerical method. For example, the predictor-corrector method above described. 
%The semi-analytical procedure for the SgKS equation is described in the following theorem.

\begin{theorem}
Consider the SgKS equation
\begin{equation}
\frac{\partial u}{\partial t}=-u \frac{\partial u}{\partial x}-\kappa \frac{\partial^2 u}{\partial x^2}-\eta \frac{\partial^3 u}{\partial x^3}-\nu\frac{\partial^4 u}{\partial x^4}+\sigma \dot{W}(t), \,\,\, x\in(a,b),\,\, t\in(0,\infty)
\label{eq:59}
\end{equation}
with constant $\sigma$, subject to the initial condition
\begin{equation}
u(x,0)=f(x), \hspace{.35cm} x\in [a,b]
\label{eq:60}
\end{equation}
and the stochastic mixed (Robin) boundary conditions 
\begin{equation}
a_1\hspace{0.05cm} u(a,t)+a_2 \hspace{0.05cm} u_x(a,t)=g_1(t),\hspace{0.2cm} t \in (0, \infty),
\label{eq:61}
\end{equation} 
\begin{equation}
b_1\hspace{0.05cm} u(b,t)+b_2 \hspace{0.05cm}u_x(b,t)=g_2(t), \hspace{0.2cm} t \in (0, \infty),
\label{eq:62}
\end{equation}
\begin{equation}
c_1\hspace{0.05cm} u(a,t)+c_2 \hspace{0.05cm} u_x(a,t)=g_3(t),\hspace{0.2cm} t \in (0, \infty),
\label{eq:63}
\end{equation} 
and
\begin{equation}
d_1\hspace{0.05cm} u(b,t)+d_2 \hspace{0.05cm}u_x(b,t)=g_4(t), \hspace{0.2cm} t \in (0, \infty),
\label{eq:64}
\end{equation}
where $a_1, a_2, b_1,c_1,c_2, d_1$ and $d_2$ are constants and $g_1(t), g_2(t), g_3(t)$ and $g_4(t)$ may be deterministic functions or stochastic functions of $t$.
Now, define
\begin {equation}
\chi=X(x,t)=x-\sigma\int\limits_{0}^{t} W(s)ds.
\label{eq:65}
\end{equation}
The solution to (\ref{eq:59})-(\ref{eq:64}) is thus given by 
\begin{equation}
u(x,t)=v\left(x-\sigma\int\limits_{0}^{t} W(s)ds,t\right)+\sigma W(t),
\label{eq:66}
\end{equation}
where $v(x,t)$ is the solution of the generalized Kuramoto--Sivashinsky equation
\begin{equation}
\frac{\partial u}{\partial t}=-u \frac{\partial u}{\partial \chi}-\kappa \frac{\partial^2 u}{\partial \chi^2}-\eta \frac{\partial^3 u}{\partial \chi^3}-\nu\frac{\partial^4 u}{\partial \chi^4}, \hspace{0.2cm} t \in (0, \infty)
\label{eq:67}
\end{equation}
for $$\chi\in\left(a-\sigma\int\limits_{0}^{t} W(s)ds,b-\sigma\int\limits_{0}^{t} W(s)ds\right)=(X(a,t),X(b,t)),$$
subject to the stochastic initial condition
\begin{equation}
v(\xi,0)=f(\xi), \hspace{0.2cm}\xi\in (X(a,t),X(b,t))
\label{eq:68}
\end{equation}
and the stochastic boundary conditions
\begin{equation}
a_1\hspace{.1cm} v(X(a,t),t)+  a_2\hspace{.1cm}v_\xi(X(a,t),t)=g_1(t)-a_1\,\sigma W(t), \hspace{0.2cm} t \in (0,\infty),
\label{eq:69}
\end{equation}
\begin{equation}
b_1\hspace{.1cm} v(X(b,t),t)+  b_2\hspace{.1cm}v_\xi(X(b,t),t)=g_2(t)-b_1\,\sigma W(t), \hspace{0.2cm} t \in (0,\infty),
\label{eq:70}
\end{equation}
\begin{equation}
c_1\hspace{.1cm} v(X(a,t),t)+  c_2\hspace{.1cm}v_\xi(X(a,t),t)=g_3(t)-c_1\,\sigma W(t), \hspace{0.2cm} t \in (0,\infty),
\label{eq:71}
\end{equation}
and 
\begin{equation}
d_1\hspace{.1cm} v(X(b,t),t)+  d_2\hspace{.1cm}v_\xi(X(b,t),t)=g_4(t)-d_1\,\sigma W(t), \hspace{0.2cm} t \in (0,\infty).
\label{eq:72}
\end{equation}
\label{thm:3}
\end{theorem}

\textbf{Proof.}  Let us define new variables $\chi=X(x,t)$ and $v(\chi,t)$ as 
\begin{equation}
\chi=X(x,t)=x-\sigma\int\limits_{0}^{t}W(s)ds
\label{eq:73}
\end{equation}
and
\begin{equation}
v(\chi,t) = u(x,t) - \sigma W(t).
\label{eq:74}
\end{equation}
Then
%%\begin{equation}
$$ \frac{\partial u}{\partial t}= \frac{\partial v}{\partial t}+ \frac{\partial v}{\partial \chi}  \frac{\partial X}{\partial t}+\sigma\dot{W}= \frac{\partial v}{\partial t}-\sigma  \frac{\partial v}{\partial \chi} W +\sigma\dot{W}$$
%% \label{bu9}
%% \end{equation}
 and
%% \begin{equation}
$$  \frac{\partial u}{\partial x}= \frac{\partial u}{\partial \chi} \frac{\partial X}{\partial x}= \frac{\partial v}{\partial \chi}, \frac{\partial^2 u}{\partial x^2}= \frac{\partial^2 v}{\partial \chi^2}, \frac{\partial^3 u}{\partial x^3}= \frac{\partial^3 v}{\partial \chi^3} \hspace{.2cm} \mbox{and} \hspace{.2cm}  \frac{\partial^4 u}{\partial x^4}= \frac{\partial^4 v}{\partial \chi^4}.$$
%% \label{bu10}
%% \end{equation}
Substituting into (\ref{eq:59}) yields
\begin{equation}
 \frac{\partial v}{\partial t}-\sigma  \frac{\partial v}{\partial \chi} W +\sigma\dot{W} =-v \frac{\partial v}{\partial \chi}-\sigma  \frac{\partial v}{\partial \chi} W-\kappa \frac{\partial^2 v}{\partial \chi^2}-\eta \frac{\partial^3 v}{\partial \chi^3}-\frac{\partial^4 v}{\partial \chi^4}+\sigma\dot{W},
\label{eq:75}
\end{equation}
which after canceling terms gives (\ref{eq:67}).
The initial condition  (\ref{eq:60}) becomes
\begin{equation}
u(x,0)=v(X(x,0),0)=v(\chi,0)=f(\chi)  
\label{eq:76}
\end{equation}
since $\chi=X(x,0)=x$, while the boundary conditions (\ref{eq:61})-(\ref{eq:64}) become, respectively,
\begin{equation}
g_1(t)=a_1\, u(a,t)+  a_2\,u_x(a,t)=a_1\, v(X(a,t),t)+  a_2\,v_\xi(X(a,t),t)+a_1\,\sigma W(t),
\label{eq:77}
\end{equation}
\begin{equation}
g_2(t)=b_1\, u(b,t)+  b_2\,u_x(b,t)=b_1\, v(X(b,t),t)+  b_2\,v_\xi(X(b,t),t)+b_1\,\sigma W(t),
\label{eq:78}
\end{equation}
\begin{equation}
g_3(t)=c_1\, u(a,t)+  c_2\,u_x(a,t)=c_1\, v(X(a,t),t)+  c_2\,v_\xi(X(a,t),t)+c_1\sigma W(t)
\label{eq:79}
\end{equation}
and 
\begin{equation}
g_4(t)=d_1\, u(b,t)+  d_2\,u_x(b,t)=d_1\, v(X(b,t),t)+  d_2\,v_\xi(X(b,t),t)+d_1\sigma W(t).
\label{eq:80}
\end{equation}
After rearranging terms, we obtain
\begin{equation}
a_1\hspace{.1cm} v(X(a,t),t)+  a_2\hspace{.1cm}v_\xi(X(a,t),t)=g_1(t)-a_1\,\sigma W(t),
\label{eq:81}
\end{equation}
\begin{equation}
b_1\hspace{.1cm} v(X(b,t),t)+  b_2\hspace{.1cm}v_\xi(X(b,t),t)=g_2(t)-b_1\,\sigma W(t), 
\label{eq:82}
\end{equation}
\begin{equation}
c_1\hspace{.1cm} v(X(a,t),t)+  c_2\hspace{.1cm}v_\xi(X(a,t),t)=g_3(t)-c_1\,\sigma W(t), 
\label{eq:83}
\end{equation}
and 
\begin{equation}
d_1\hspace{.1cm} v(X(b,t),t)+  d_2\hspace{.1cm}v_\xi(X(b,t),t)=g_4(t)-d_1\,\sigma W(t)
\label{eq:84}
\end{equation}
which are (\ref{eq:69}) and (\ref{eq:72}), respectively.
\qed

Thus, we numerically solve the IBVP (\ref{eq:67})-(\ref{eq:72}) involving a deterministic PDE using an appropriate classical numerical method rather than solving the IBVP (\ref{eq:59})-(\ref{eq:64}) which involves a stochastic PDE. Having numerically computed $v(\chi,t)$, (\ref{eq:66}) is then used to obtain $u(x,t)=v(\chi(x,t),t)+\sigma W(t)$. 

The main issue is that in a numerical computation, two independent time variables have to be taken into consideration in the implementations of the initial and boundary conditions of $v$ , one in the expression for $\chi=X(x,t)$ and another in the expression for $v(\chi,t)$. In addition, a very small time step is needed in the numerical simulations. Thus, the implementation of the semi-analytical procedure is computationally expensive and requires a very large allocation of computer memory. 

\section{WCE based numerical solutions of the SKS and SgKS equations}\label{sec:5}

In this section, we apply the WCE based numerical method described in section \ref{subsec:4.1} to numerically solve some IBVPs involving SKS equation (\ref{eq:1}) and the SgKS equation (\ref{eq:2}).
%As mentioned before in section \ref{sec:1} (Introduction), the SKS equation is obtained by setting $\eta=0$ in the SgKS equation (\ref{eq:33}). 
In the numerical simulations, $\sigma$ is a constant, set to $\sigma=1$, in order the semi-analytical solution to be useful. The results presented in the present paper were obtained using the time step $\Delta t=0.005$ and the spatial step $\Delta x=0.2$ in the numerical computations.  

Our truncation is quite simple and contains Gaussian terms only. The Wick polynomials are such that $T_\alpha=H_{1}(\xi_i), i=\alpha, \alpha\le\tilde{I}$,  $T_\alpha=H_{\alpha_i=\alpha}(\xi_i), i=\alpha, \tilde{I}<\alpha\le I$, $\sum_i\alpha_i<I$, where $I$ is the number of terms kept in the Paley-Wiener representation of $W(t)$ (\ref{eq:10}) in the numerical implementation, and represents the order of the WCE.  This corresponds to the WCE  truncation,
\begin{equation}
u_{W,I}(x,t)=\sum\limits_{\alpha\le I}u_\alpha (x,t) T_\alpha= \sum\limits_{\alpha=\alpha_i=i\le \tilde{I}}u_{\alpha} (x,t)\xi_i\\+\sum\limits_{\tilde{I}\le\alpha=\alpha_i=i\le I}u_\alpha (x,t)H_{\alpha_i}(\xi_i),
\label{eq:85}
\end{equation}
and the double sum in the propagator (\ref{eq:44}) has to be truncated accordingly.  This truncation is simple, Gaussian and captures important stochastic information provided in the Paley-Wiener representation of $W(t)$. In the WCE approximate numerical solution $\tilde{I}=40$ and $I=60$.
 
To assess the accuracy of the WCE approximate numerical solution, the numerical solution is contrasted with a semi-analytical solution obtained following the procedure described in section \ref{subsec:4.2}.
%in which the  problem is transformed by a change of independent variables into an equivalent problem  comprising the deterministic Burgers' equation with stochastic boundary and initial conditions which is then solved numerically.
In  each test problem, the absolute difference between the numerical solution $u_{W,I}$ and the semi-analytical solution is calculated according to (\ref{eq:30}) for each $x=x_k = a +k \Delta x$  ($k = 0, 1, \dots, K$) and $t=t_n = n \Delta t$ ($n=  0, 1, \dots, N$) and  their relative difference  over the interval $x \in [a,b]$ is calculated  according to (\ref{eq:32}) for each $t=t_n = n \Delta t$. In the formulas (\ref{eq:30}) and (\ref{eq:32}), the semi-analytical solution is used in place of the analytical solution $u$. The  WCE numerical computations took 942.5 sec CPU time compare with 24.0 sec CPU time for the semi-analytical solution on a Pentium (R) PC with 2.5 GHz CPU even though the semi-analytical solution needs a very large allocation of computer memory.

\subsection{Numerical solution of the SKS equation}\label{subsec:5.1}

We examine two test problems. In each test problem, an IBVPs involving the SKS equation (\ref{eq:2}) is numerically solved. The constant parameters $\kappa, \eta$ and $\nu$ are respectively set to $\kappa=0.1, \eta=0$ and $\nu=0.02$.

\subsubsection{Test problem 1:}\label{prob:1} The SKS equation (\ref{eq:2}) is solved numerically on the domain $[a,b]\times[0,T]=[-10,10]\times[0,3]$ subject to the initial condition $u(x,0) = f(x)=\frac{\cos{\left({\pi x}/20\right)}}{3.5+\sin{\left({\pi x}/20\right)}}$, the stochastic boundary conditions $u(-10,t)=u(10,t)=\sigma W(t)$; and for simplification purpose, the other two boundary conditions are assumed to be periodic. This stochastic boundary condition is implemented because the stochastic properties of the solutions have to be taken into account at the boundaries of the domain. With the spatial step size $\Delta x=0.2$ and the time step  size $\Delta t= 0.005$, the space interval is divided into 100 points and the time interval into 600 points. 

The associate propagator has the initial condition
\begin{equation}
u_\alpha(x,0)=\begin{cases}\frac{\cos{\left(\frac{\pi x}{20}\right)}}{3.5+\sin{\left(\frac{\pi x}{20}\right)}}, & \text{if $\alpha=0$,} \\
0, & \text{if $\alpha\ne 0$, } \\
\end{cases}\\
\label{eq:86}
\end{equation}
and the boundary conditions 
\begin{equation}
u_{\alpha}(-10,t)=u_{\alpha}(10,t)=\begin{cases}
0, & \text{if $\alpha=0$}, \\
\sigma\sum\limits_{i=1}^{\infty}\mathbb{I}_{\alpha_j=\delta_{i,j}} \int_0^t m_i(s) ds, & \text{if $\alpha\ne 0$}. \\
\end{cases}\\
\label{eq:87}
\end{equation}

\begin{figure}[!htb]
\begin{picture}(0,0)
\put (-10,180){(a)}
\put (-25,95){$u(x,t)$}
\put (120,-5){$x$}
\put (-10,-25){(b)}
\put (-25,-115){$u(x,t)$}
\put (120,-205){$x$}
\put (-10,-225){(c)}
\put (-25,-305){$u(x,t)$}
\put (120,-405){$x$}
\end{picture}
\vspace{.75cm}
\centering{\includegraphics[scale=0.60]{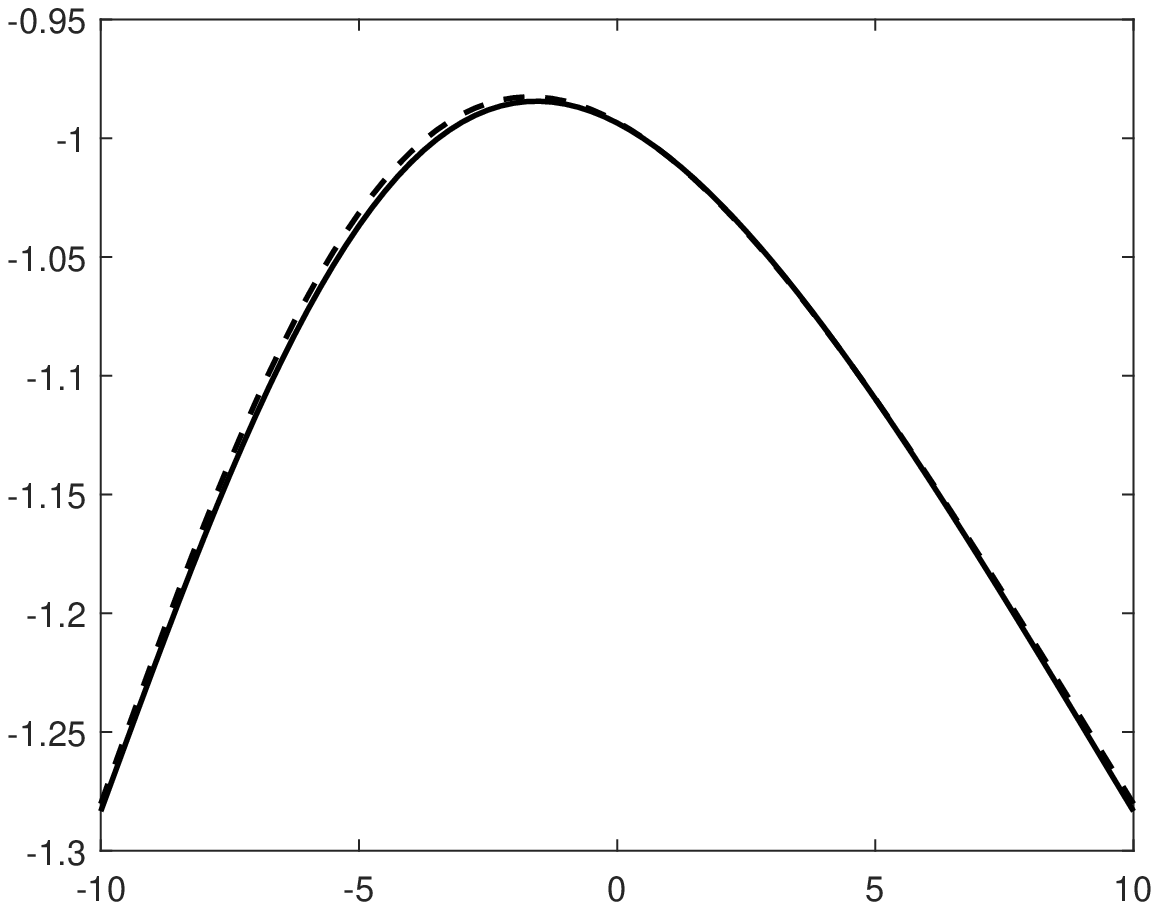}} \vspace{.75cm}
\centering{\includegraphics[scale=0.60]{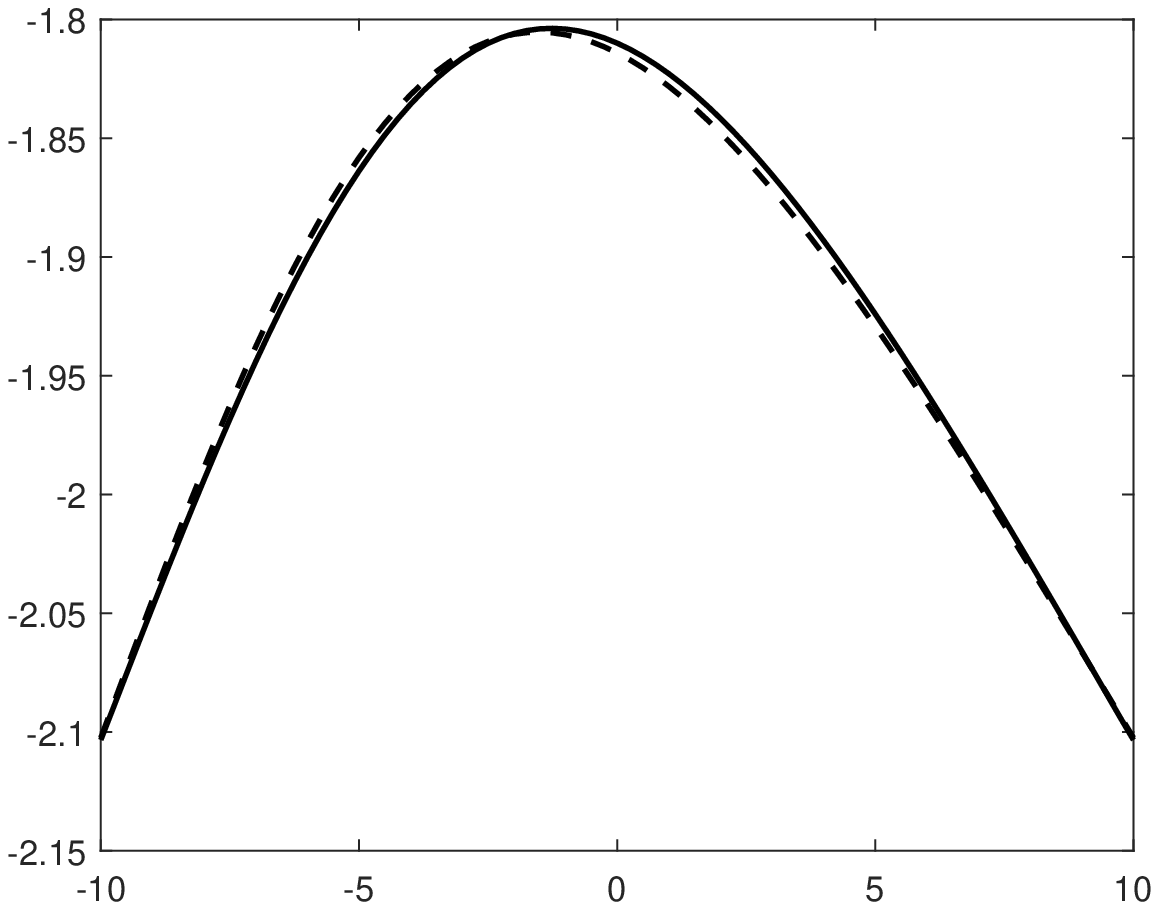}} \vspace{.25cm}
\centering{\includegraphics[scale=0.60]{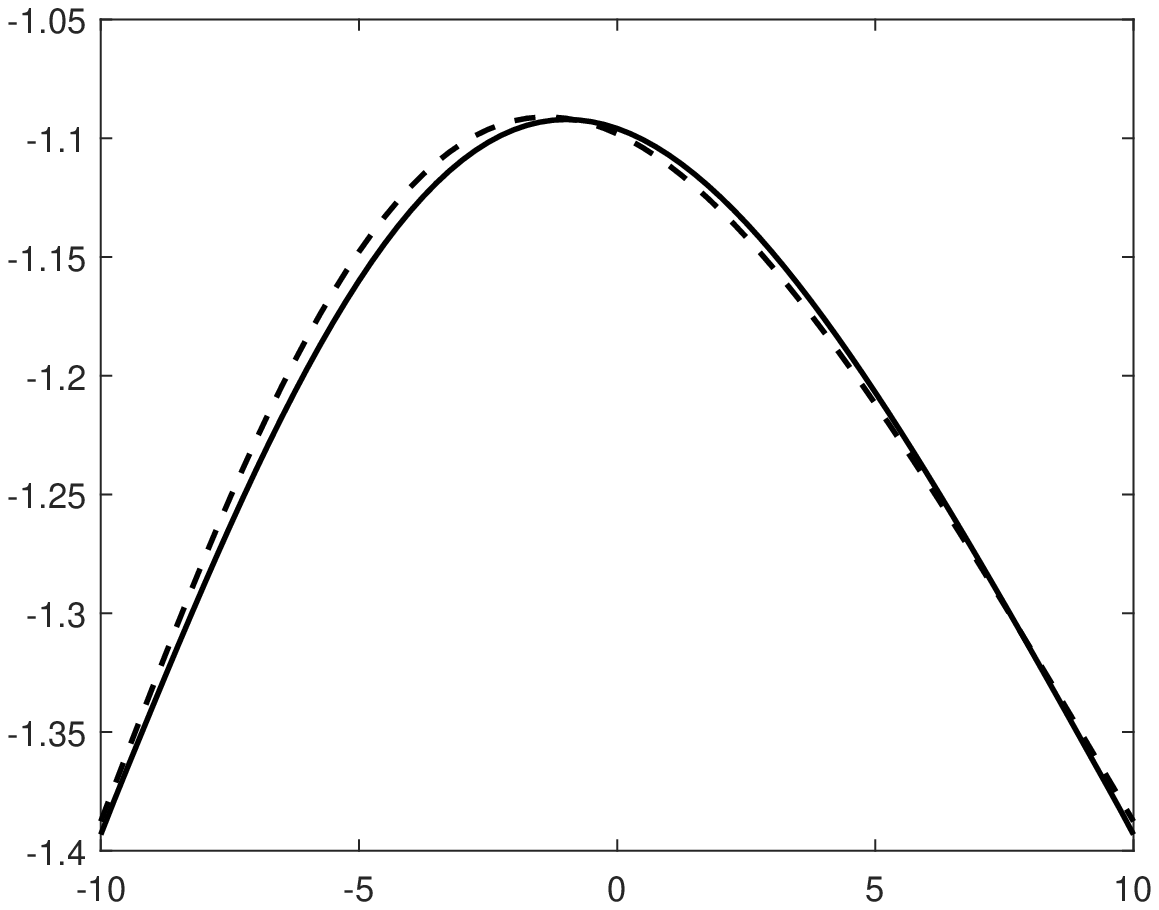}}
\caption{Stochastic Kuramoto--Sivashinsky (SKS) equation (\ref{eq:2}): solution $u(x,t)$ as a function of position $x$ at the time, (a) $t=1$, (b) $t=2$ and (d) $t=3$. The solid line is the semi-analytical solution while the dashed line is the WCE-based numerical solution.}
\label{fig:5}
\end{figure}

\begin{figure}[!htb]
\begin{picture}(0,0)
\put (-25,90){$u(x,t)$}
\put (125,-5){$t$}
\end{picture}
\centering{\includegraphics[scale=0.6]{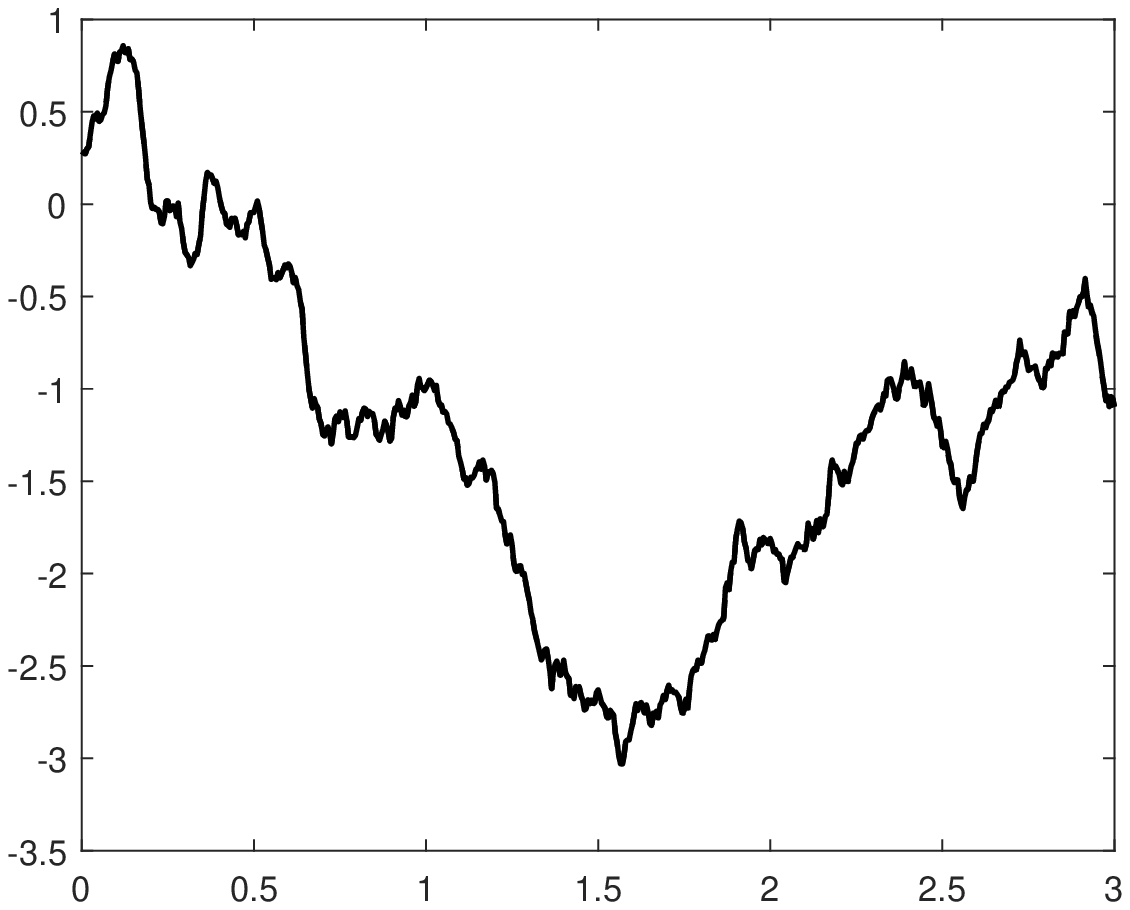}}
\caption{Stochastic Kuramoto--Sivashinsky (SKS) equation (\ref{eq:2}): solution $u(x,t)$ as a function of time $t$ at position $x=0$. The solid line is the semi-analytical solution while the dashed line is the WCE-based numerical solution.}
\label{fig:6}
\end{figure}

\begin{figure}[!htb]
\begin{picture}(0,0)
\put (5,150){(a)}
\put (-20,75){$\Delta_a u$}
\put (100,-5){$t$}
\put (220,150){(b)}
\put (200,75){$\Delta_r u$}
\put (318,-5){$t$}
\end{picture}
\centering
{\includegraphics[scale=0.5]{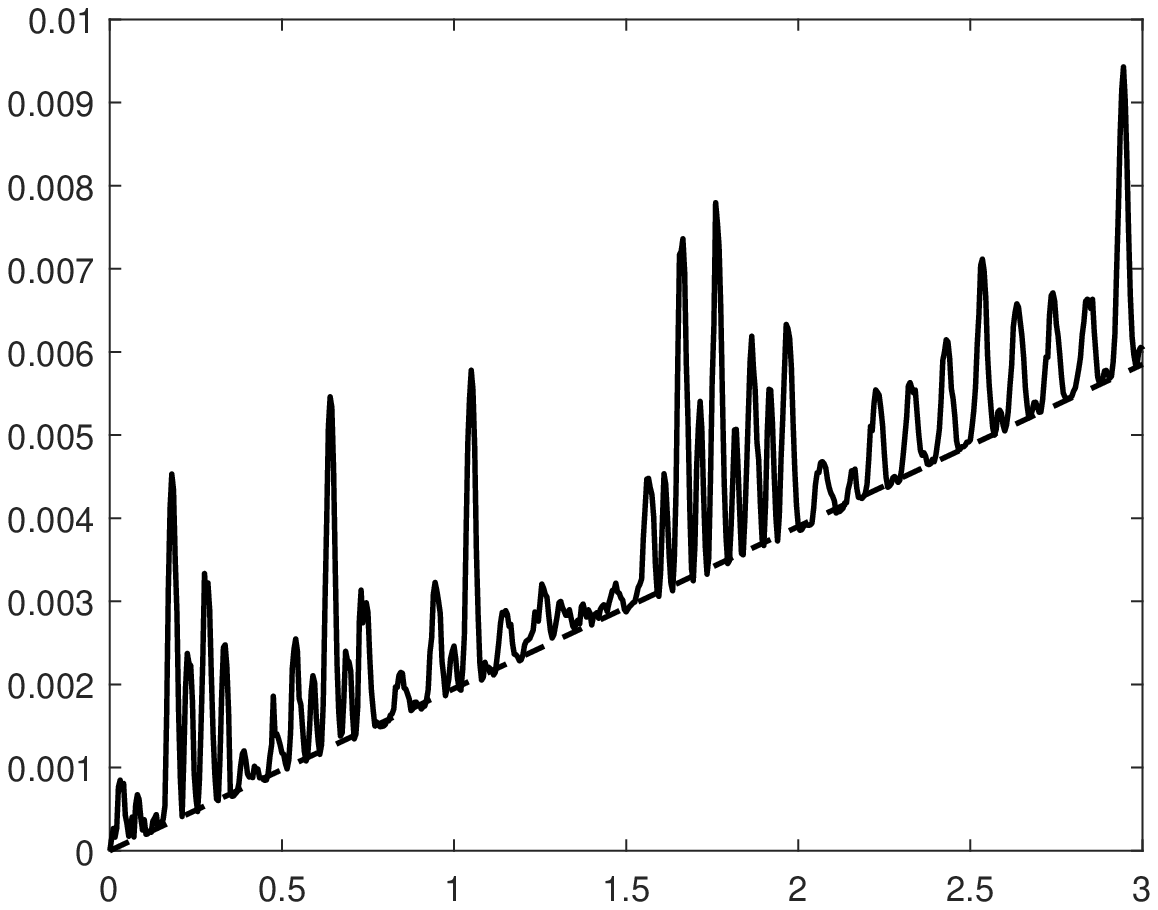}\hspace{.8cm} \includegraphics[scale=0.5]{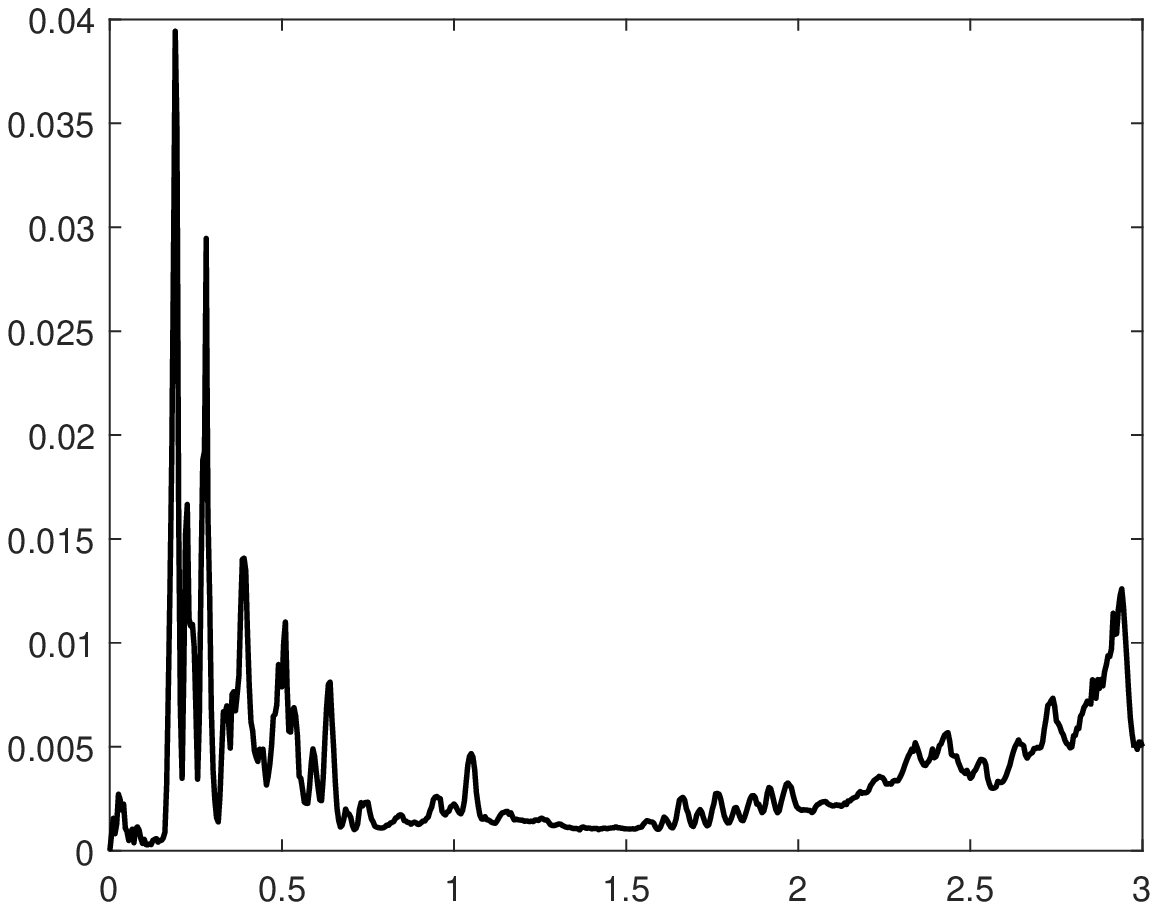}}
\caption{Stochastic Kuramoto--Sivashinsky (SKS) equation (\ref{eq:2}): (a) absolute error $\Delta_a u$ as a function of time $t$, and (b) relative error $\Delta_r u$ as a function of time $t$.}
\label{fig:7}
\end{figure}

Some results are shown in Figures \ref{fig:5}-\ref{fig:7}. It is seen in these figures that the semi-analytical solution and the WCE approximate numerical solution are in good agreement and their absolute difference $\Delta_a u$ is order $10^{-3}$ or less for $t$ in the interval $[0,3]$. It is seen in Figure \ref{fig:7}(a) that the absolute error $\Delta_a u\sim O(\varsigma t)$, where $\varsigma \sim O(10^{-3})$. Therefore, the absolute error should be $O(10^{-2})$ by $t=10$. It is shown in Figure \ref{fig:7}(b) that the relative (error) $\Delta_r u$ is $O(10^{-2})$ or less  for all $t$ in the interval $[0,3]$ and its maximum value on this interval is $4\%$. 

\subsubsection{Test problem 2:} \label{prob:2} The SKS equation (\ref{eq:2}) is solved numerically on the domain $[a,b]\times[0,T]=[0,20]\times[0,3]$ subject to the initial condition $f(x)=\frac{\sin{\left({\pi x}/20\right)}-\sin{\left({\pi x}/10\right)}}{7.5-\cos{\left({\pi x}/20\right)}+0.5\cos{\left({\pi x}/10\right)}}$ and the stochastic boundary conditions $u(0,t)=u(20,t)=\sigma W(t)$. For simplification purpose, the other two boundary conditions are assumed to be periodic as in the test problem 1.

The initial condition for the propagator is 
\begin{equation}
u_\alpha(x,0)=\begin{cases}
\frac{\sin{\left(\frac{\pi x}{20}\right)}-\sin{\left(\frac{\pi x}{10}\right)}}{7.5-\cos{\left(\frac{\pi x}{20}\right)}+0.5\cos{\left(\frac{\pi x}{10}\right)}}, & \text{if $\alpha=0$,} \\
0, & \text{if $\alpha\ne 0$, } \\
\end{cases}\\
\label{eq:88}
\end{equation}
and the boundary conditions are
\begin{equation}
u_{\alpha}(0,t)=u_{\alpha}(20,t)=\begin{cases}
0, & \text{if $\alpha=0$}, \\
\sigma\sum\limits_{i=1}^{\infty}\mathbb{I}_{\alpha_j=\delta_{i,j}}  \int_0^t m_i(s) ds, & \text{if $\alpha\ne 0$}. \\
\end{cases}\\
\label{eq:89}
\end{equation}

\begin{figure}[!htb]
\begin{picture}(0,0)
\put (-10,180){(a)}
\put (-25,95){$u(x,t)$}
\put (120,-5){$x$}
\put (-10,-25){(b)}
\put (-25,-105){$u(x,t)$}
\put (120,-205){$x$}
\put (-10,-225){(c)}
\put (-25,-305){$u(x,t)$}
\put (120,-405){$x$}
\end{picture}
\centering{\includegraphics[scale=0.6]{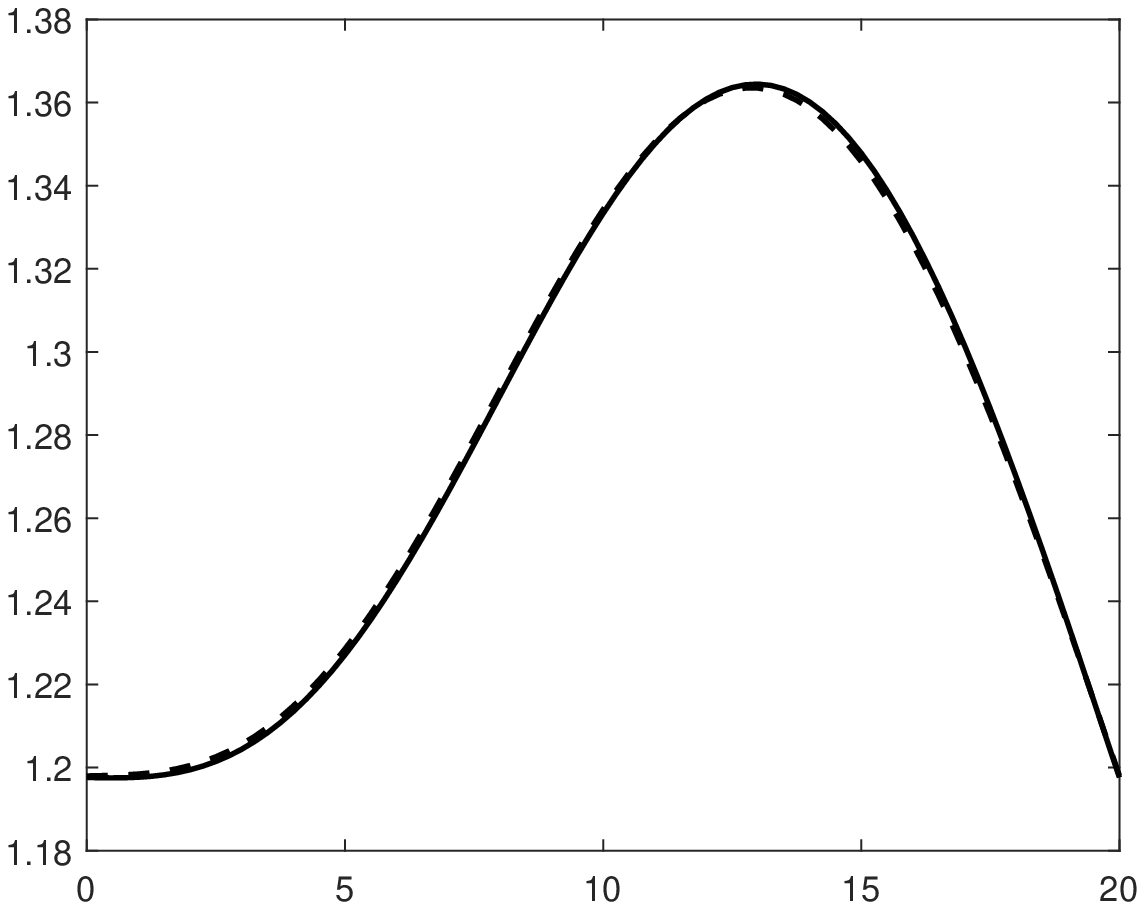}}\vspace{.75cm}
\centering{\includegraphics[scale=0.6]{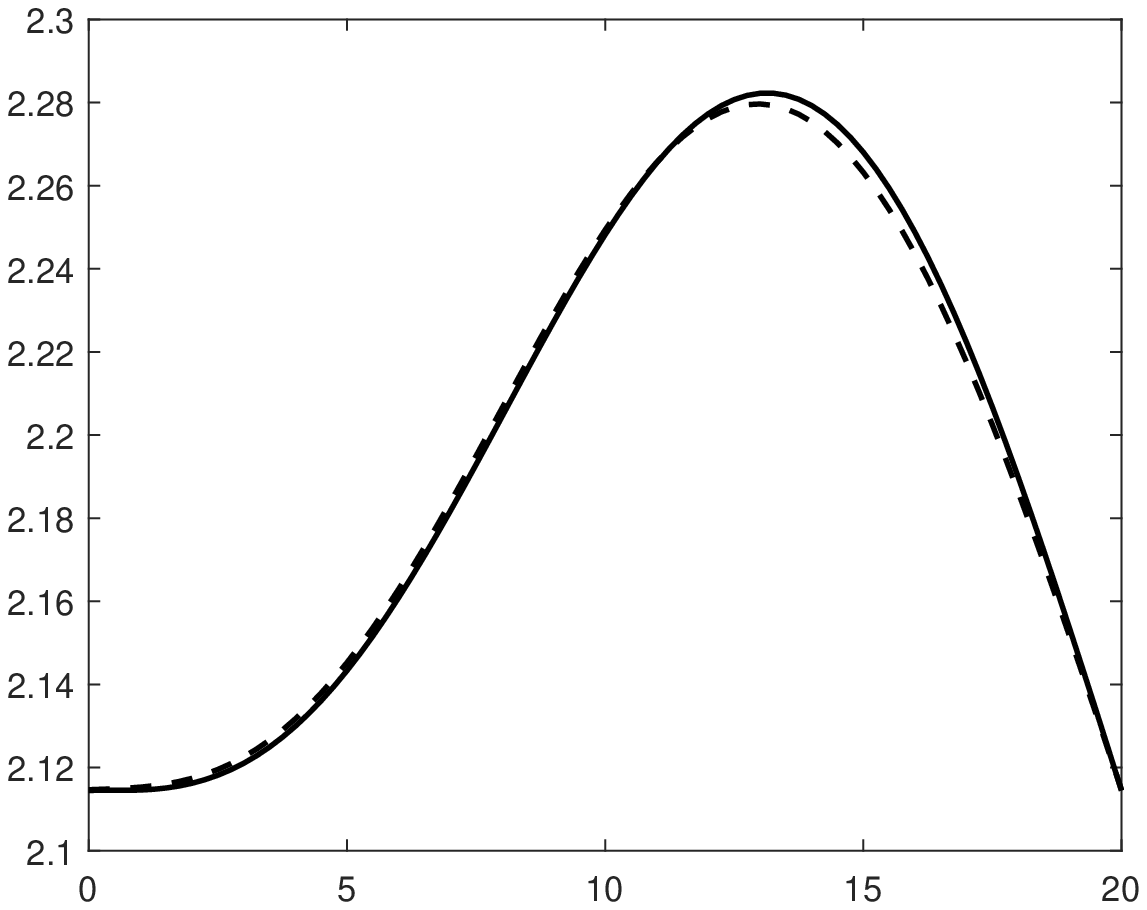}}\vspace{.75cm}
\centering{\includegraphics[scale=0.6]{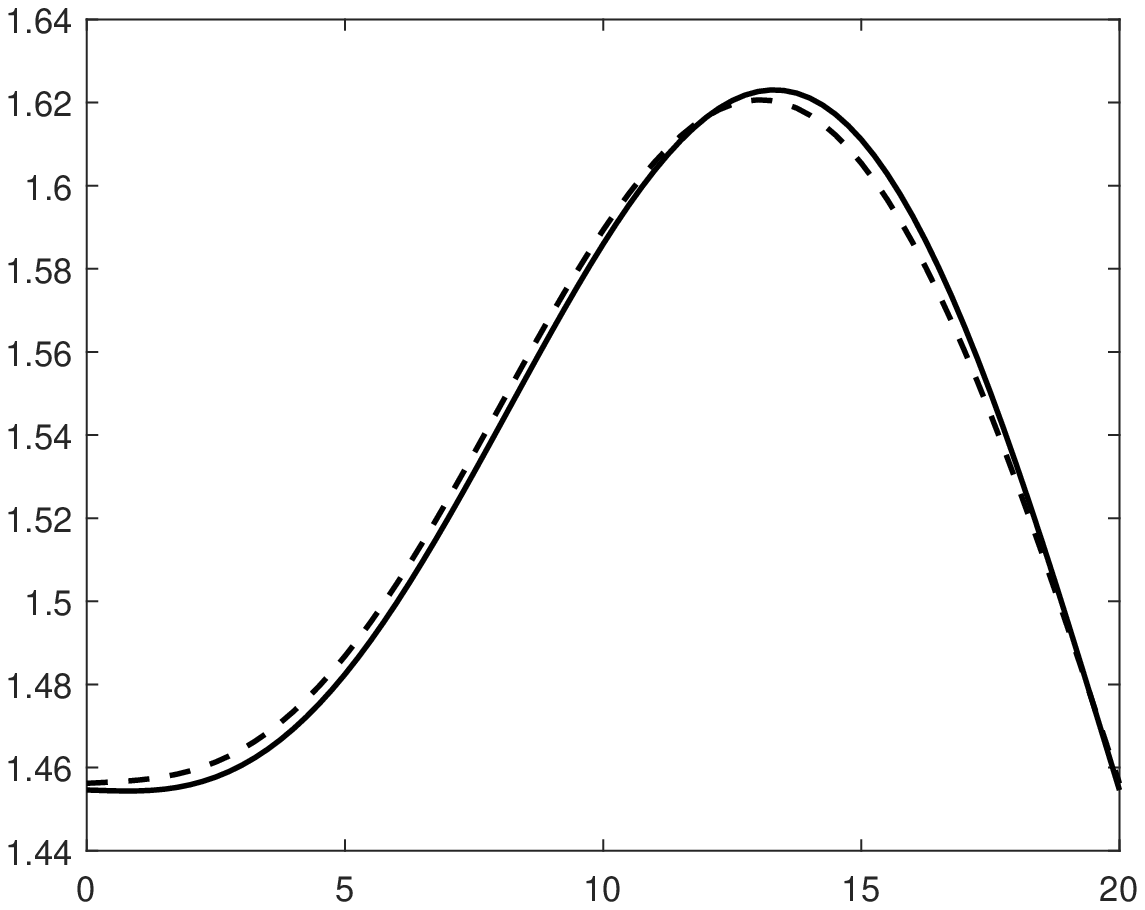}}
\caption{Stochastic Kuramoto--Sivashinsky (SKS) equation (\ref{eq:2}): solution $u(x,t)$ as a function of position $x$ at the time, (a) $t=1$, (b) $t=2$ and (d) $t=3$. The solid line is the semi-analytical solution while the dashed line is the WCE-based numerical solution.}
\label{fig:8}
\end{figure}

\begin{figure}[!htb]
\begin{picture}(0,0)
\put (-25,90){$u(x,t)$}
\put (125,-5){$t$}
\end{picture}
\centering{\includegraphics[scale=0.6]{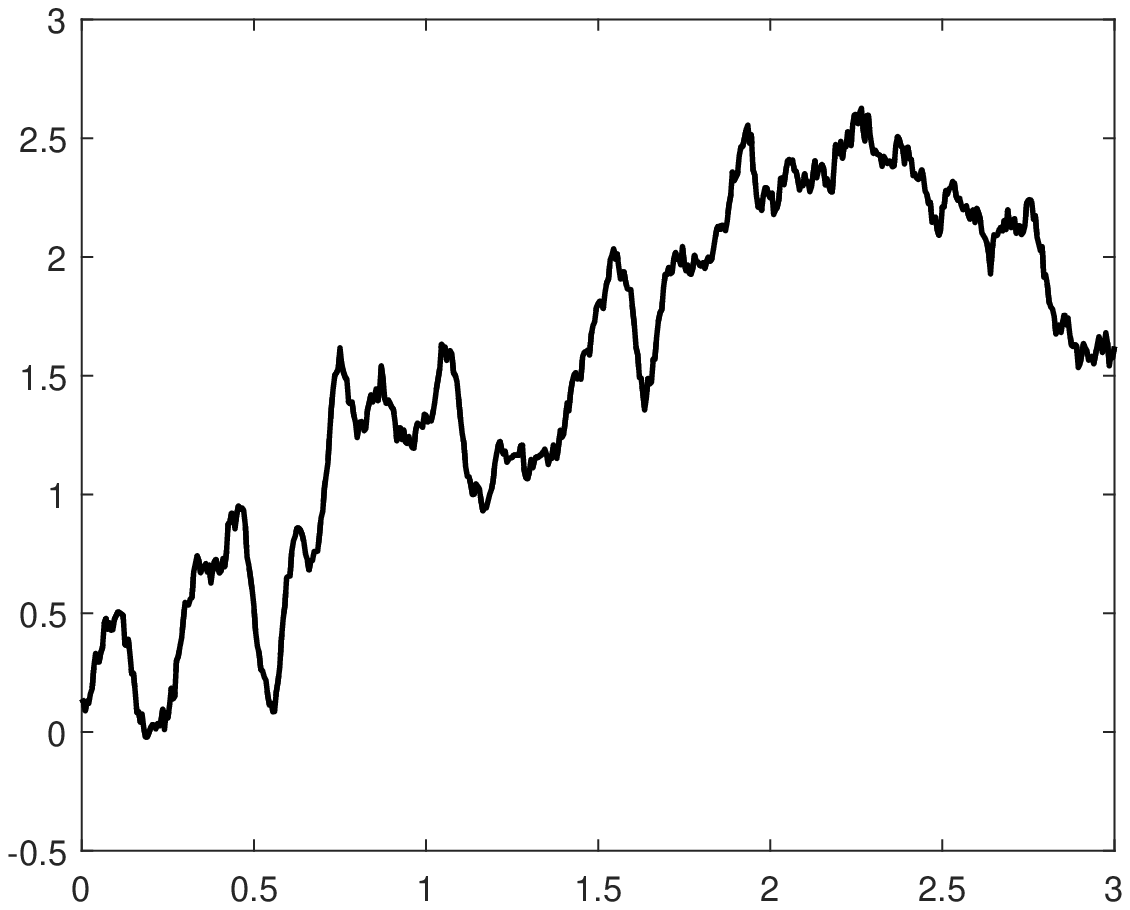}}
\caption{Stochastic Kuramoto--Sivashinsky (SKS) equation (\ref{eq:2}):  solution $u(x,t)$ as a function of time $t$ at position $x=10$. The solid line is the semi-analytical solution while the dashed line is the WCE-based numerical solution.}
\label{fig:9}
\end{figure}

\begin{figure}[!htb]
\begin{picture}(0,0)
\put (5,150){(a)}
\put (-20,75){$\Delta_a u$}
\put (100,-5){$t$}
\put (220,150){(b)}
\put (200,75){$\Delta_r u$}
\put (318,-5){$t$}
\end{picture}
\centering{\includegraphics[scale=0.5]{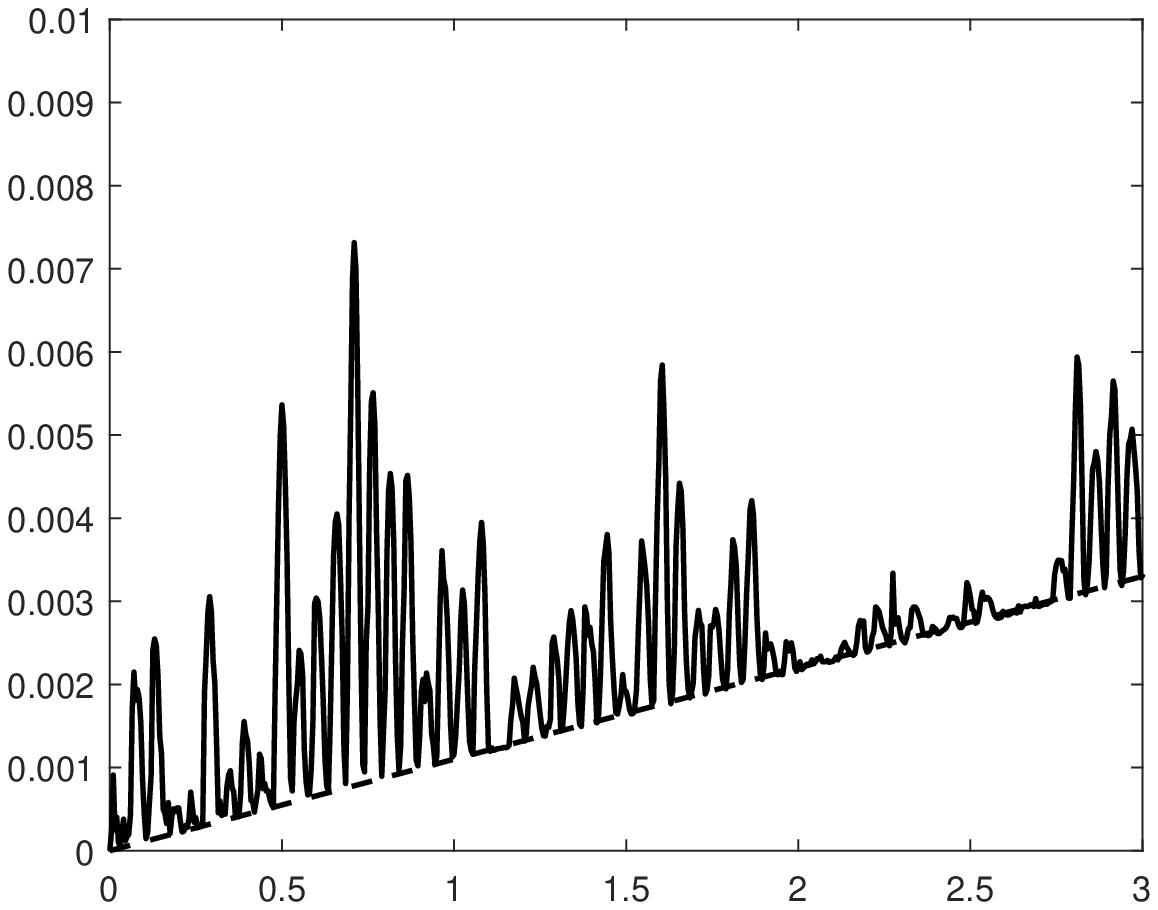}\hspace{.8cm}\includegraphics[scale=0.5]{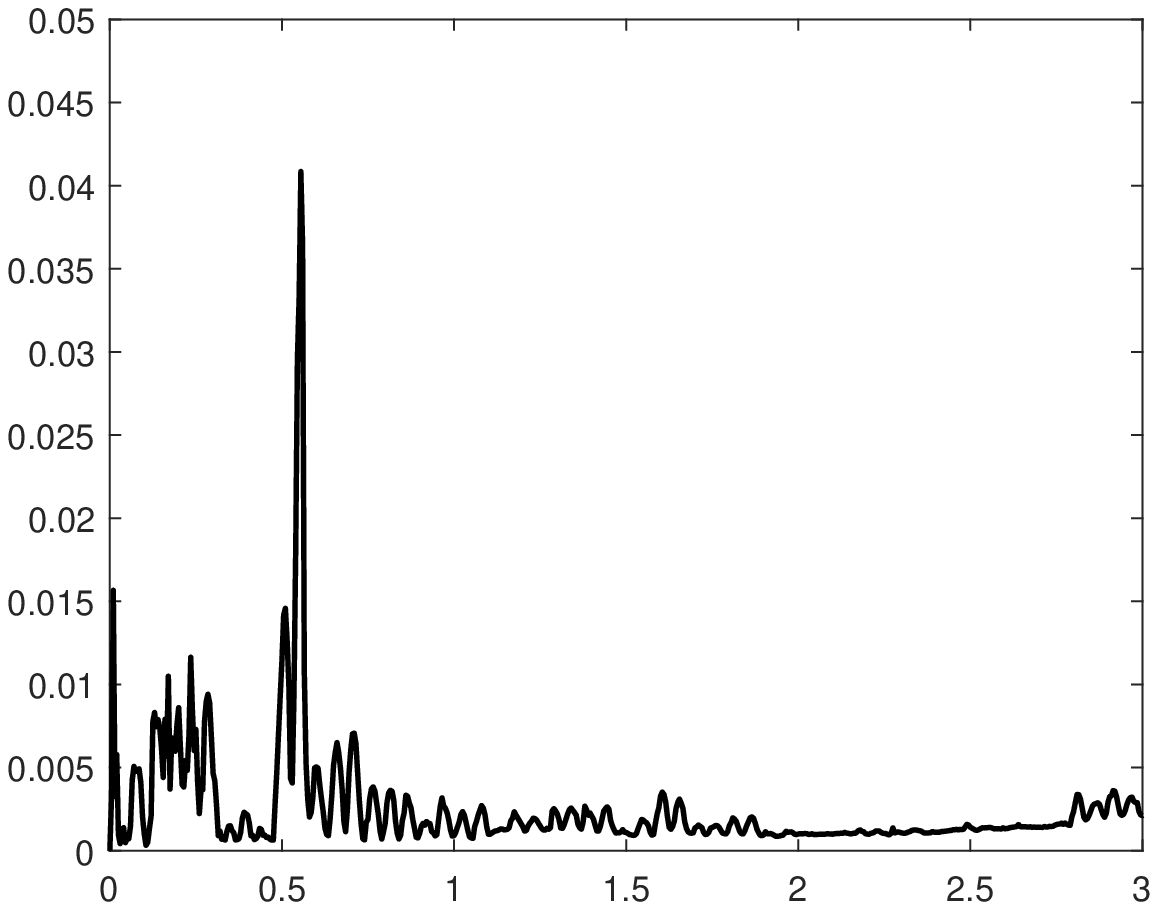}}
\caption{Stochastic Kuramoto--Sivashinsky (SKS) equation (\ref{eq:2}): (a) absolute error $\Delta_a u$ as a function of time $t$, and (b) relative error $\Delta_r u$ as a function of time $t$.}
\label{fig:10}
\end{figure}

Some results are shown in Figures \ref{fig:8}-\ref{fig:10}. These figures show that the semi-analytical solution and the WCE approximate numerical solution are in good agreement and their absolute difference $\Delta_a u$ is order $10^{-3}$ or less for $t$ in the interval $[0,3]$. Figure \ref{fig:10}(a) shows that the absolute error $\Delta_a u\sim O(\varsigma t)$, where $\varsigma \sim O(10^{-3})$, and so, the absolute error should be $O(10^{-2})$ by $t=10$. It is seen in Figure \ref{fig:10}(b) that the relative (error) $\Delta_r u$ is $O(10^{-2})$ or less  for all $t$ in the interval $[0,3]$ and its maximum value on this interval is $4\%$.

\subsection{Numerical solution of the SgKS equation, $\eta\ne0$}\label{subsec:5.1}

Two test problems are also examined. In each test problem, an IBVP involving the SgKS equation (\ref{eq:1}) is numerically solved. The constant parameters $\kappa, \eta$ and $\nu$ are respectively set to $\kappa=0.1, \eta=0.05$ and $\nu=0.02$. The values of the parameters $\kappa$ and $\nu$ are the same as in the test problems 1 and 2 except that $\eta$ is now nonzero.

\subsubsection{Test problem 3:}\label{prob:3} The SgKS equation (\ref{eq:2}) is numerically solved on the domain $[a,b]\times[0,T]=[-10,10]\times[0,3]$. The initial and boundary conditions are similar to those in test problem 1, $u(x,0) = f(x)=\frac{\cos{\left({\pi x}/20\right)}}{3.5+\sin{\left({\pi x}/20\right)}}$, $u(-10,t)=u(10,t)=\sigma W(t)$; and for simplification purpose, the other two boundary conditions are assumed to be periodic. In that case, the initial condition for the propagator (\ref{eq:44}) is given by (\ref{eq:86}), while the boundary condition is given (\ref{eq:87}).

\begin{figure}[!htb]
\begin{picture}(0,0)
\put (-10,180){(a)}
\put (-25,95){$u(x,t)$}
\put (120,-5){$x$}
\put (-10,-25){(b)}
\put (-25,-115){$u(x,t)$}
\put (120,-205){$x$}
\put (-10,-225){(c)}
\put (-25,-305){$u(x,t)$}
\put (120,-405){$x$}
\end{picture}
%\vspace{.75cm}
\centering{\includegraphics[scale=0.6]{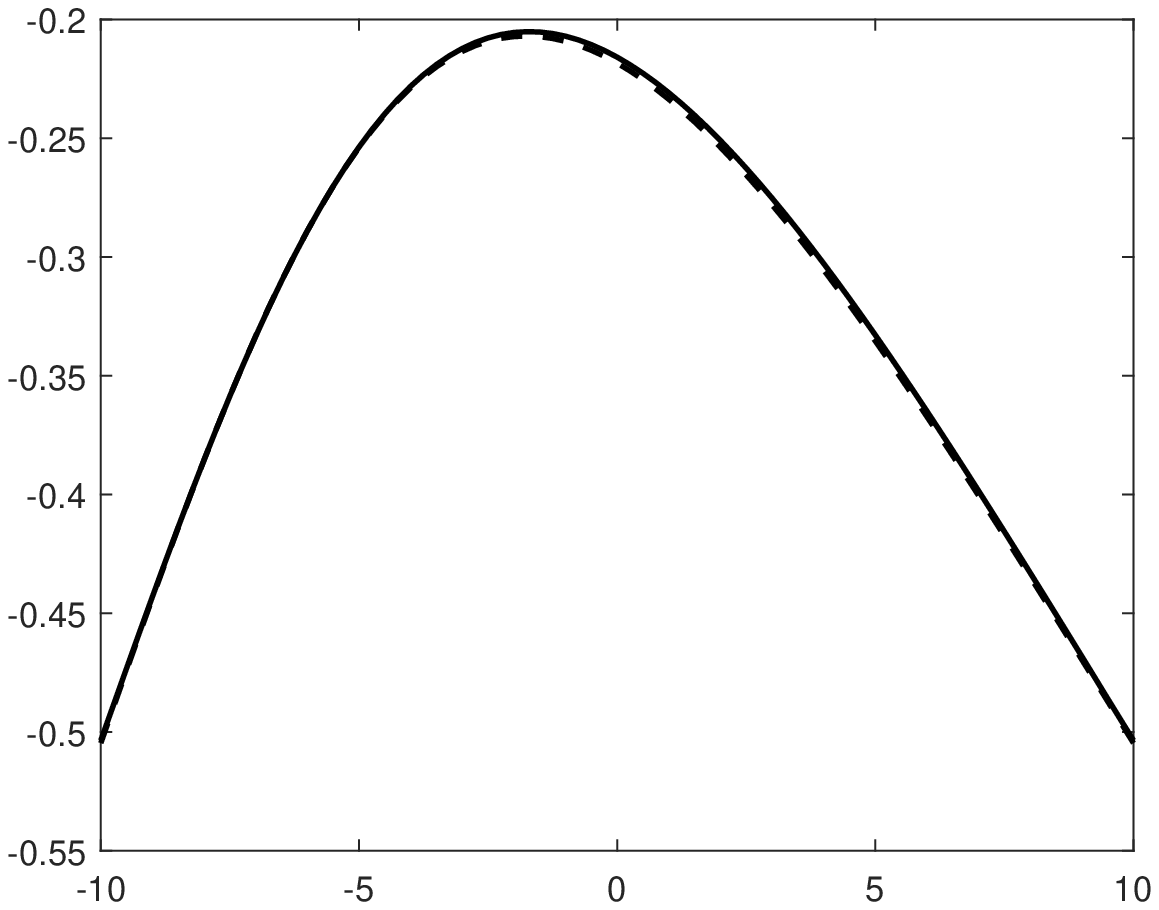}}\vspace{.75cm}
\centering{\includegraphics[scale=0.6]{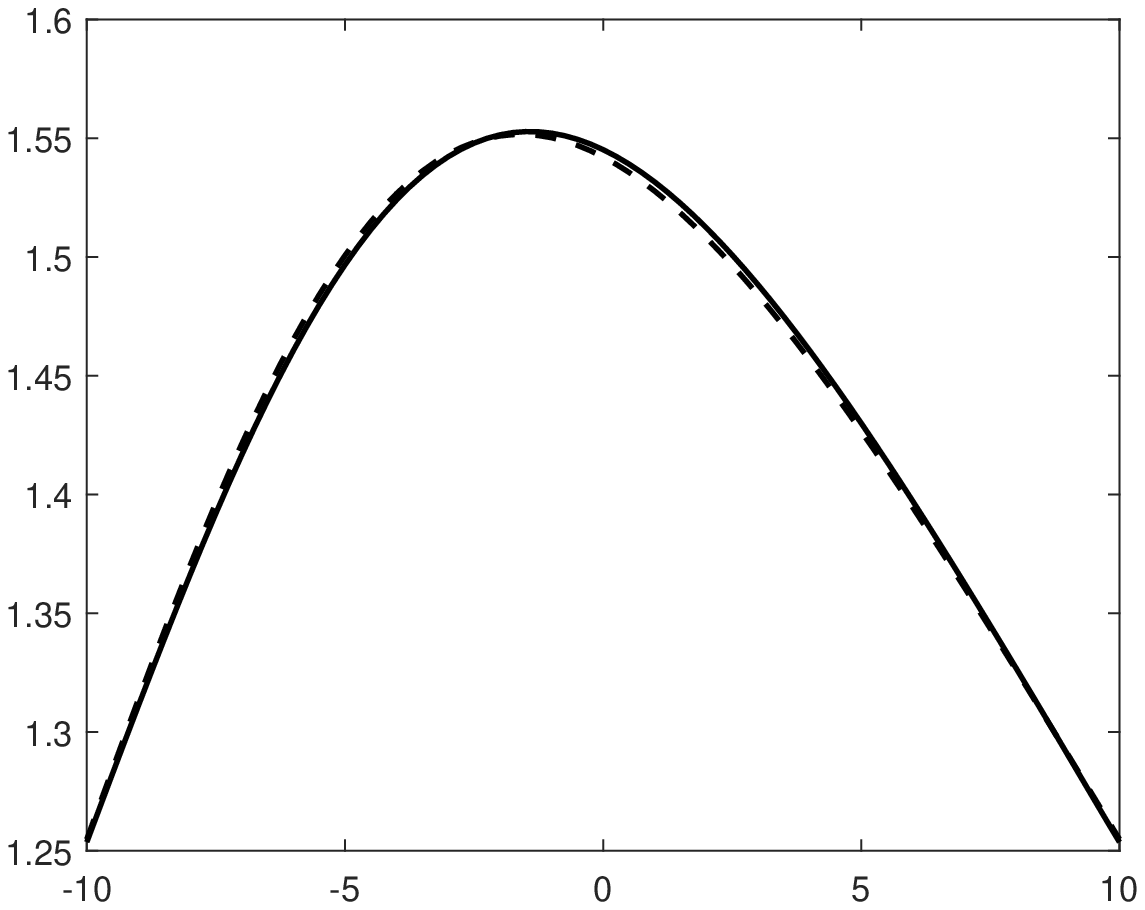}}\vspace{.75cm}
\centering{\includegraphics[scale=0.6]{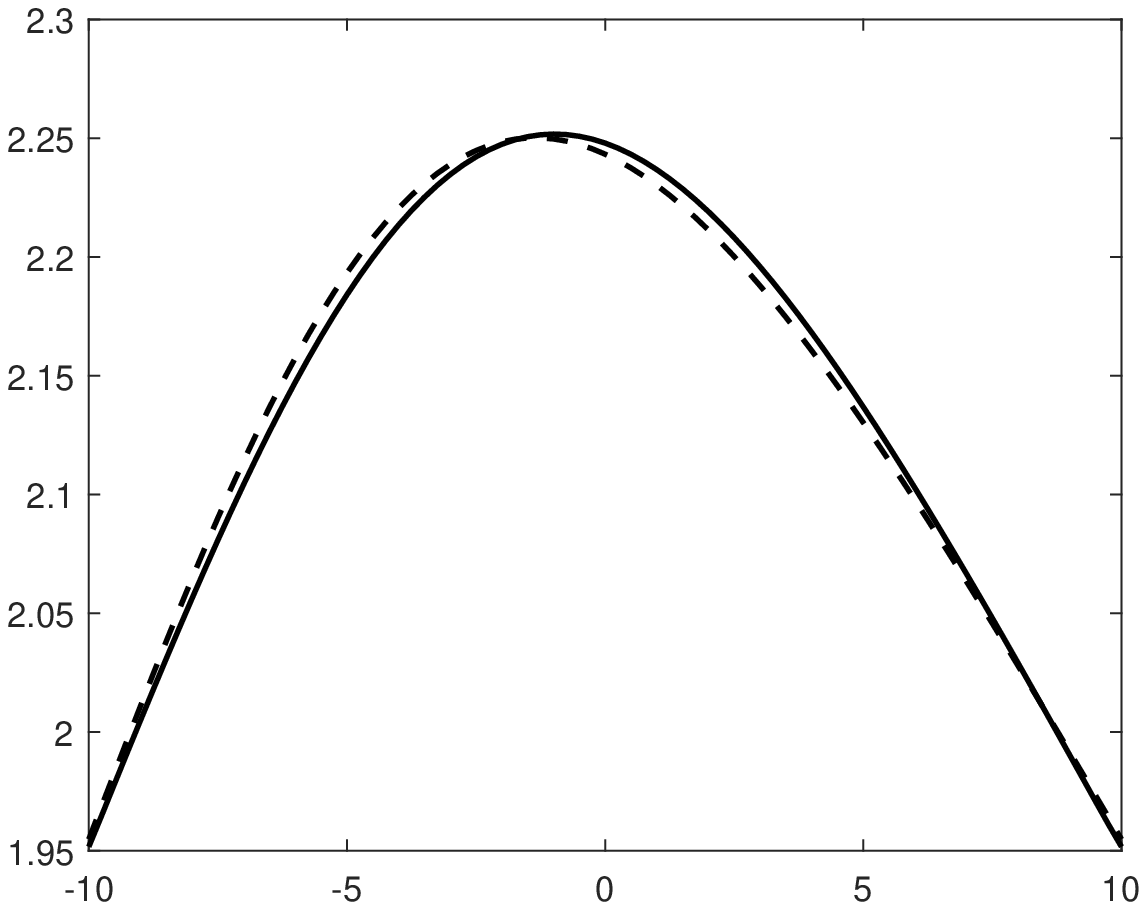}}
\caption{Stochastic generalized  Kuramoto--Sivashinsky (SgKS) equation (\ref{eq:1}): solution $u(x,t)$ as a function of position $x$ at the time, (a) $t=1$, (b) $t=2$ and (d) $t=3$. The solid line is the semi-analytical solution while the dashed line is the WCE-based numerical solution.}
\label{fig:11}
\end{figure}

\begin{figure}[!htb]
\begin{picture}(0,0)
\put (-25,90){$u(x,t)$}
\put (125,-5){$t$}
\end{picture}
\centering{\includegraphics[scale=0.6]{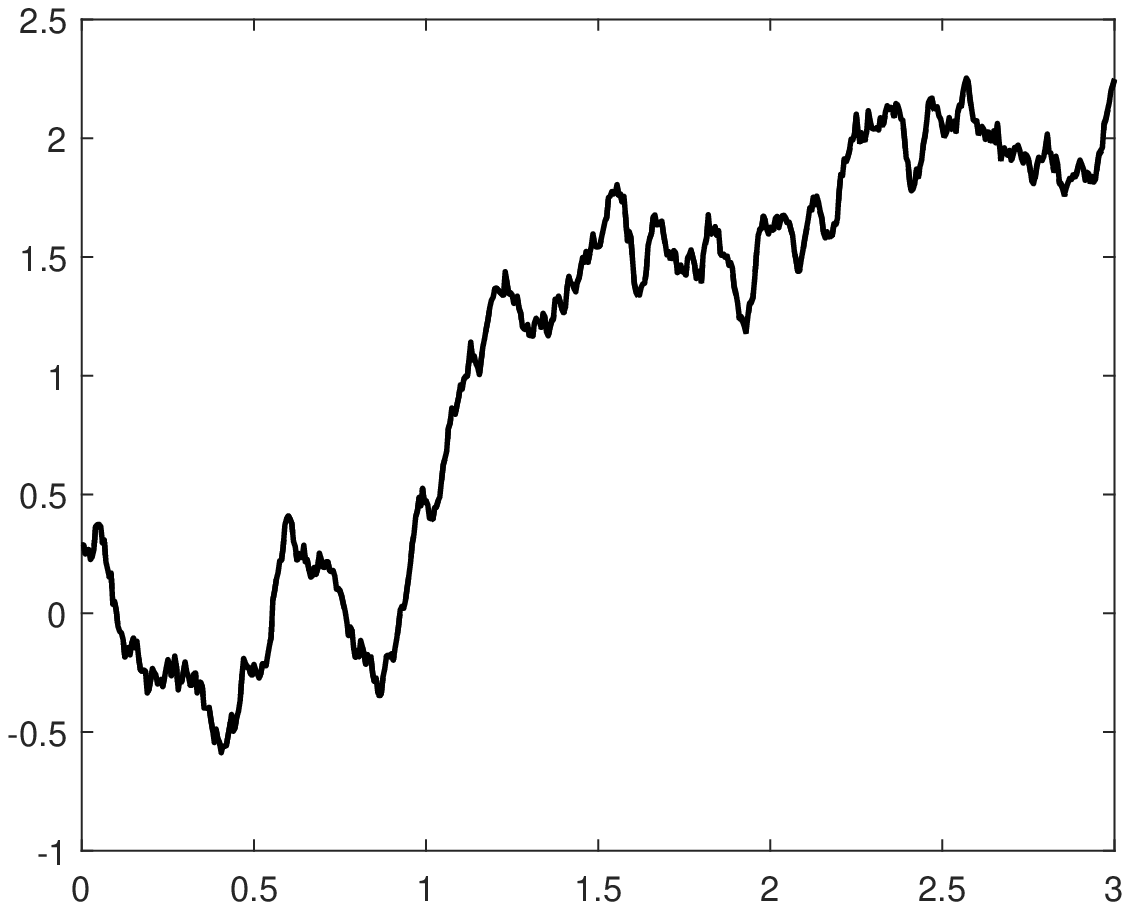}}
\caption{Stochastic generalized  Kuramoto--Sivashinsky (SgKS) equation (\ref{eq:1}):  solution $u(x,t)$ as a function of time $t$ at position $x=0$. The solid line is the semi-analytical solution while the dashed line is the WCE-based numerical solution.}
\label{fig:12}
\end{figure}

\begin{figure}[!htb]
\begin{picture}(0,0)
\put (5,150){(a)}
\put (-20,75){$\Delta_a u$}
\put (100,-5){$t$}
\put (220,150){(b)}
\put (200,75){$\Delta_r u$}
\put (318,-5){$t$}
\end{picture}
\centering{\includegraphics[scale=0.5]{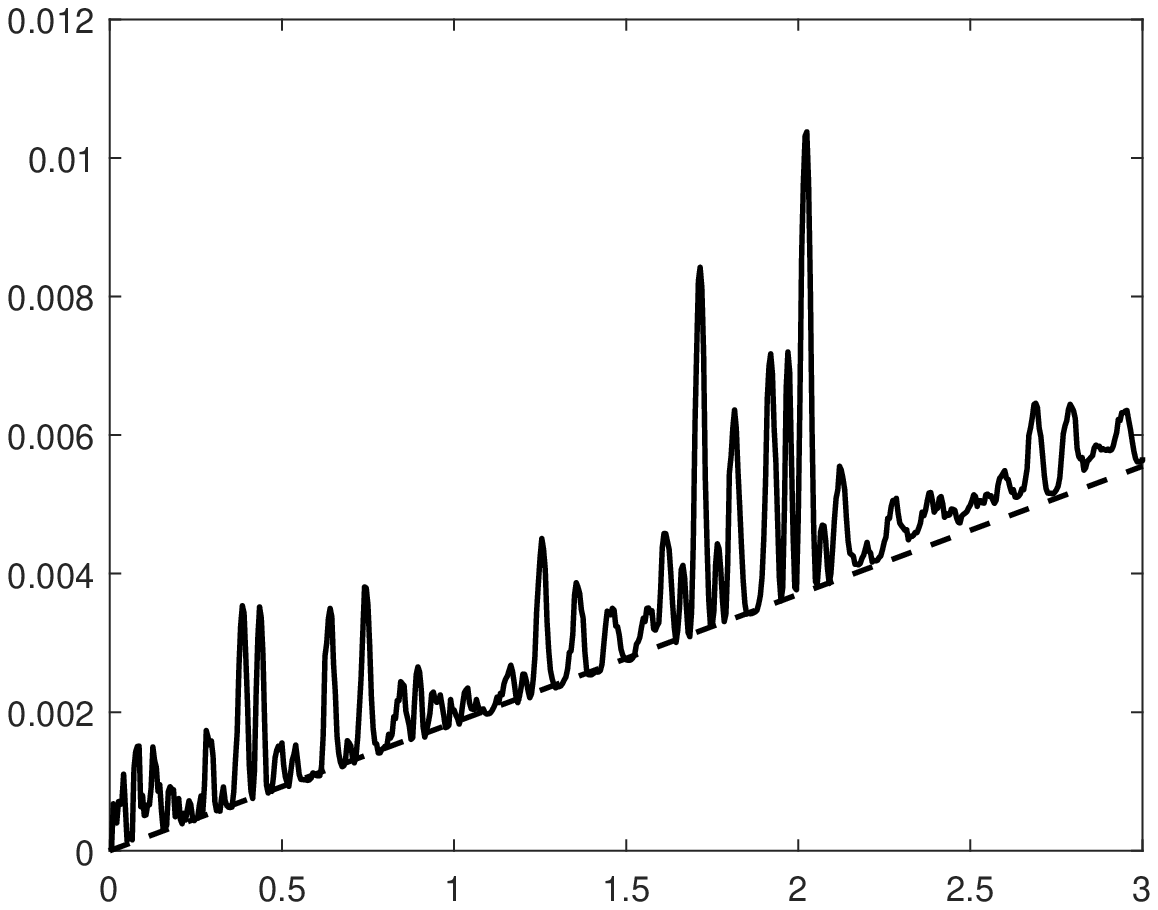}\hspace{.8cm}\includegraphics[scale=0.5]{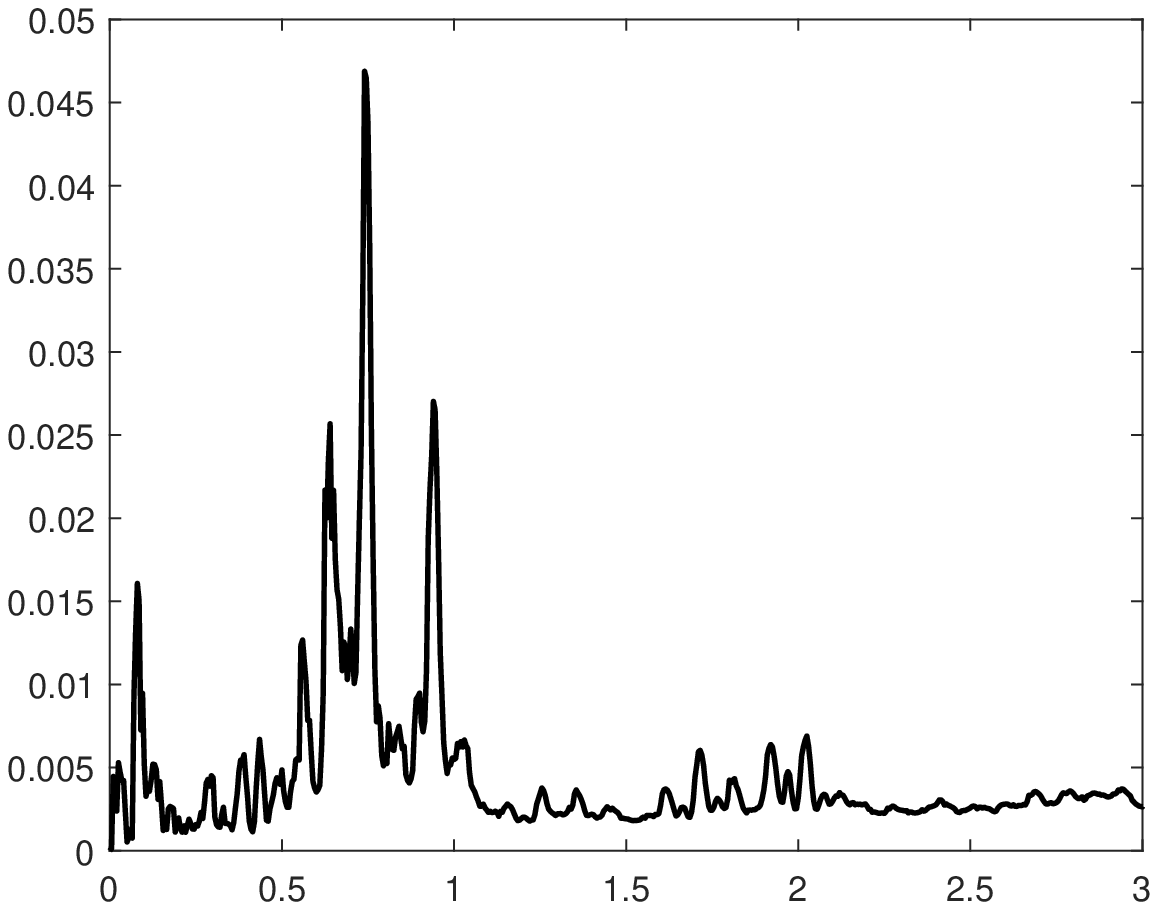}}
\caption{Stochastic generalized  Kuramoto--Sivashinsky (SgKS) equation (\ref{eq:1}): (a) absolute error $\Delta_a u$ as a function of time $t$, and (b) relative error $\Delta_r u$ as a function of time $t$.}
\label{fig:13}
\end{figure}

Some results are presented in Figures \ref{fig:11}-\ref{fig:13}. These figures show that the semi-analytical solution and the WCE approximate numerical solution are in good agreement. Although, their absolute difference (error) $\Delta_a u$ may reach a value of 0.01, it is order $10^{-3}$ in general on the time interval $[0,3]$. As seen in Figure \ref{fig:13}(a), the absolute error follows a linear pattern and is $O(\varsigma t)$ as in the test problems 1 and 2, where $\varsigma \sim O(10^{-3})$. The relative (error) is shown in Figure \ref{fig:13}(b). It is order $10^{-2}$ or less for all $t$ in the interval $[0,3]$ as in test problems 1 and 2, and its maximum value on this interval is $4.5\%$. 

\subsubsection{Test problem 4:} \label{prob:4} The SgKS equation (\ref{eq:2}) is  numerically solved on the domain $[a,b]\times[0,T]=[0,20]\times[0,3]$ subject to the initial condition $f(x)=\frac{\sin{\left({\pi x}/20\right)}-\sin{\left({\pi x}/10\right)}}{7.5-\cos{\left({\pi x}/20\right)}+0.5\cos{\left({\pi x}/10\right)}}$ and the stochastic boundary conditions are $u(0,t)=u(20,t)=\sigma W(t)$ as in the test problem 2; and for simplification purpose, the other two boundary conditions are also assumed to be periodic as in the test problem 2. Therefore, the initial condition for the propagator (\ref{eq:44}) is given by (\ref{eq:88}), while the boundary condition is given (\ref{eq:89}) .

\begin{figure}[!htb]
\begin{picture}(0,0)
\put (-10,180){(a)}
\put (-25,95){$u(x,t)$}
\put (120,-5){$x$}
\put (-10,-25){(b)}
\put (-25,-115){$u(x,t)$}
\put (120,-205){$x$}
\put (-10,-225){(c)}
\put (-25,-305){$u(x,t)$}
\put (120,-405){$x$}
\end{picture}
%\vspace{.75cm}
\centering{\includegraphics[scale=0.6]{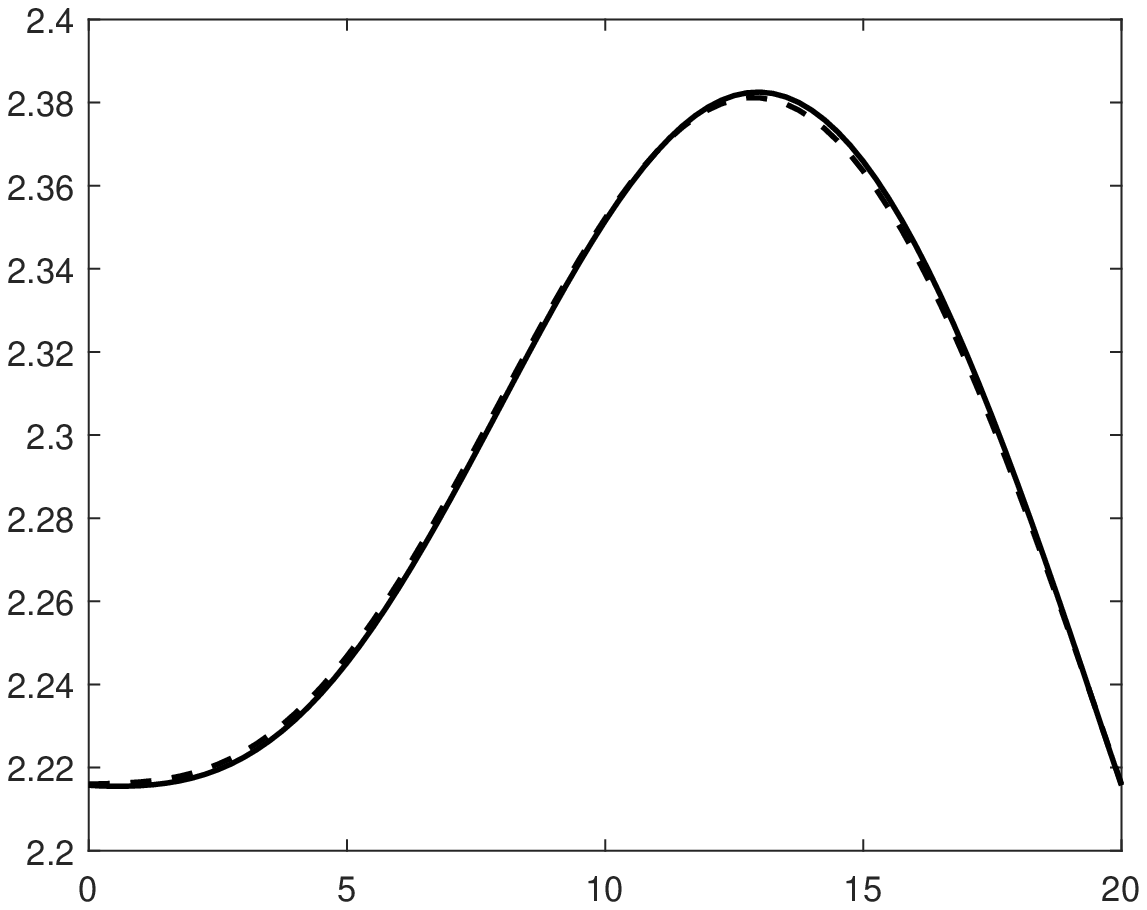}}\vspace{.75cm}
\centering{\includegraphics[scale=0.6]{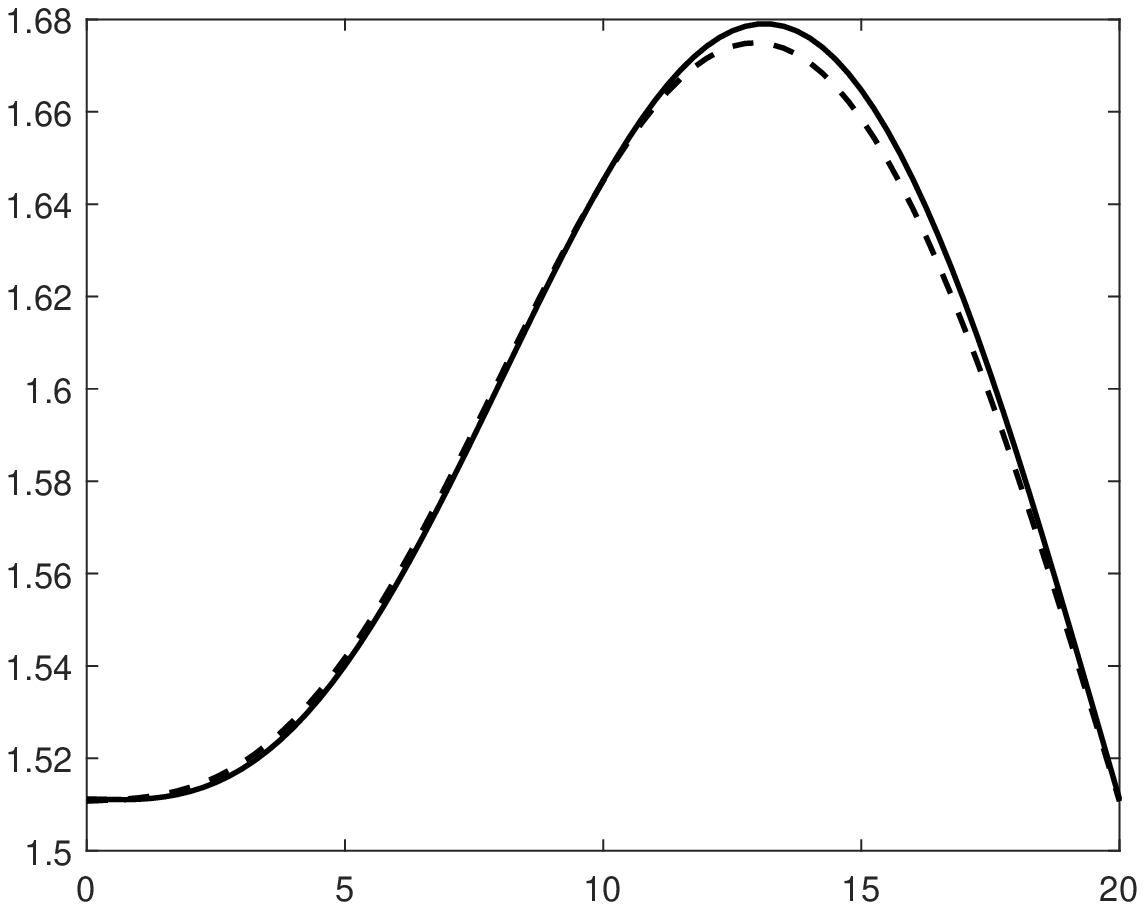}}\vspace{.75cm}
\centering{\includegraphics[scale=0.6]{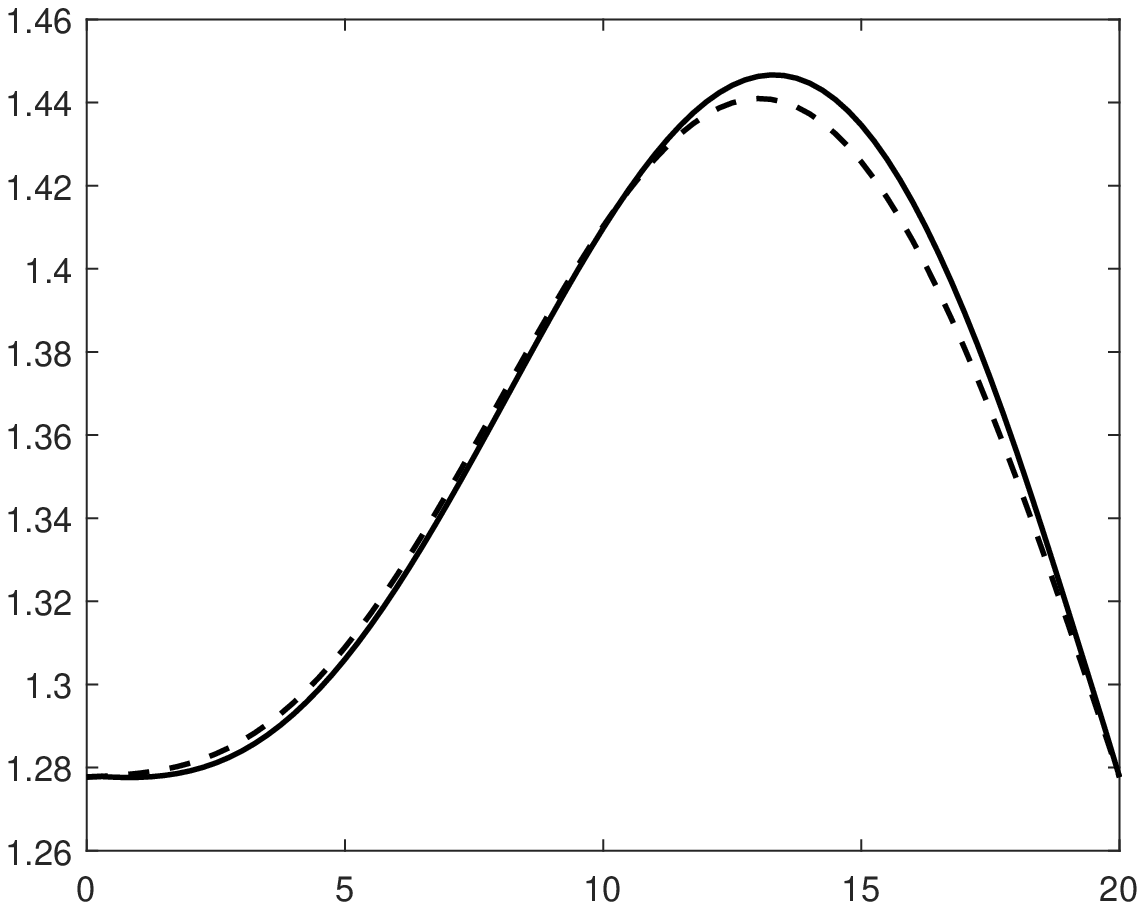}}
\caption{Stochastic generalized  Kuramoto--Sivashinsky (SgKS) equation (\ref{eq:1}): solution $u(x,t)$ as a function of position $x$ at the time, (a) $t=1$, (b) $t=2$ and (d) $t=3$. The solid line is the semi-analytical solution while the dashed line is the WCE-based numerical solution.}
\label{fig:14}
\end{figure}

\begin{figure}[!htb]
\begin{picture}(0,0)
\put (-25,90){$u(x,t)$}
\put (125,-5){$t$}
\end{picture}
\centering{\includegraphics[scale=0.6]{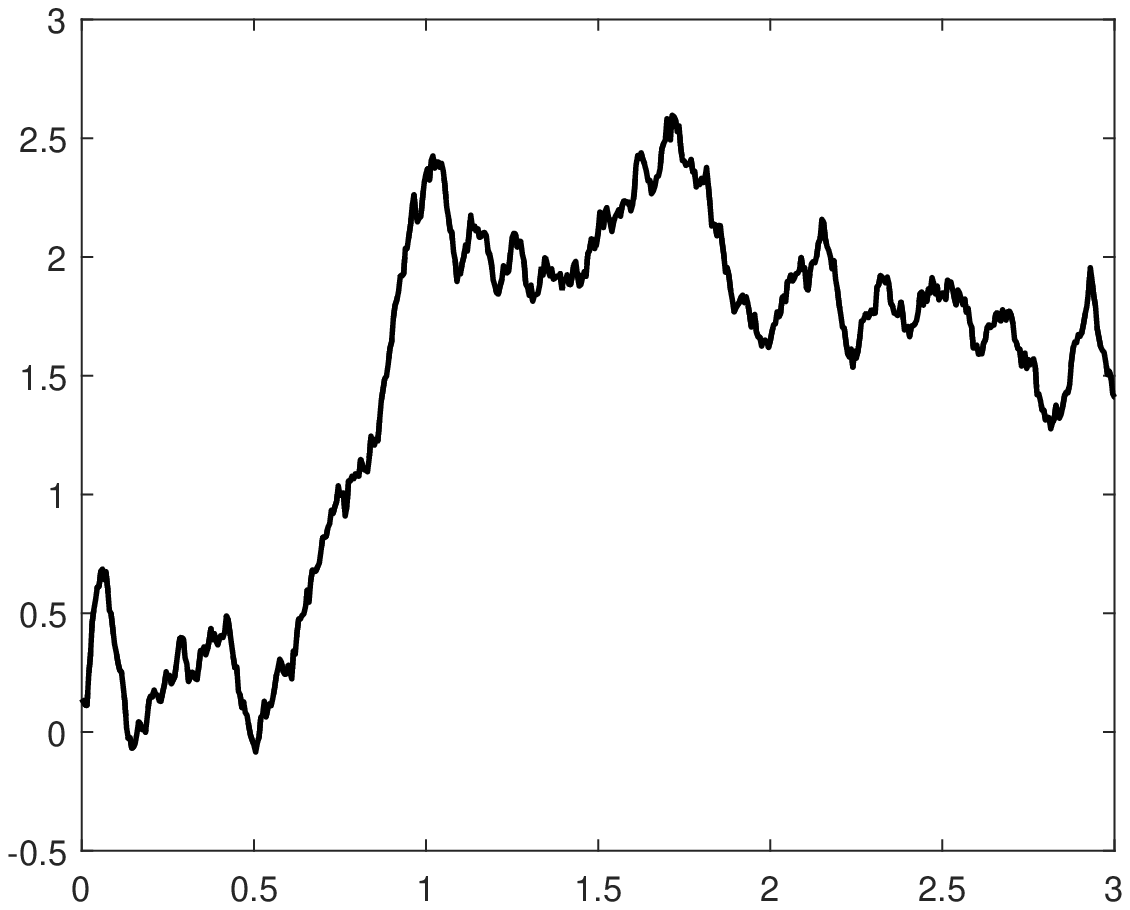}}
\caption{Stochastic generalized  Kuramoto--Sivashinsky (SgKS) equation (\ref{eq:1}): solution $u(x,t)$ as a function of time $t$ at position $x=10$. The solid line is the semi-analytical solution while the dashed line is the WCE-based numerical solution.}
\label{fig:15}
\end{figure}

\begin{figure}[!htb]
\begin{picture}(0,0)
\put (5,150){(a)}
\put (-20,75){$\Delta_a u$}
\put (100,-5){$t$}
\put (220,150){(b)}
\put (200,75){$\Delta_r u$}
\put (318,-5){$t$}
\end{picture}
\centering{\includegraphics[scale=0.5]{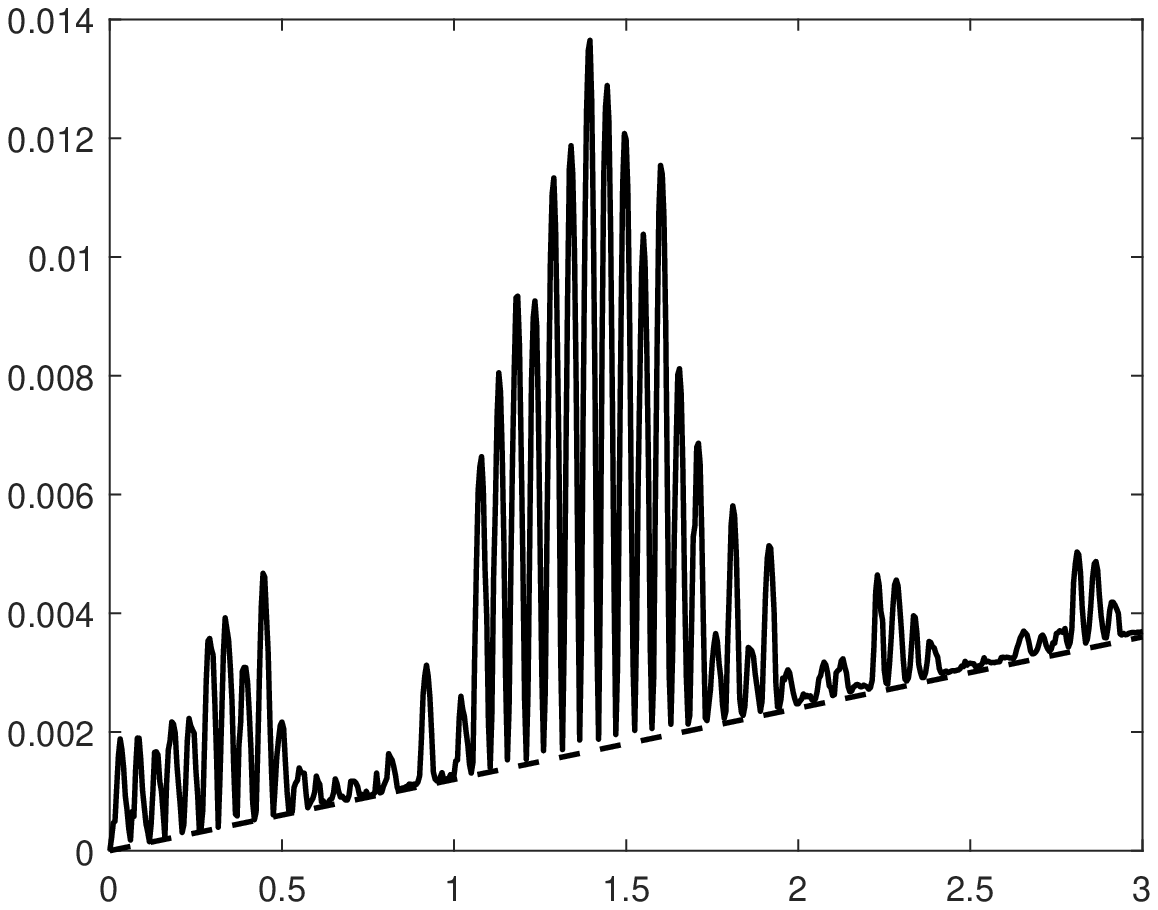}\hspace{.8cm}\includegraphics[scale=0.5]{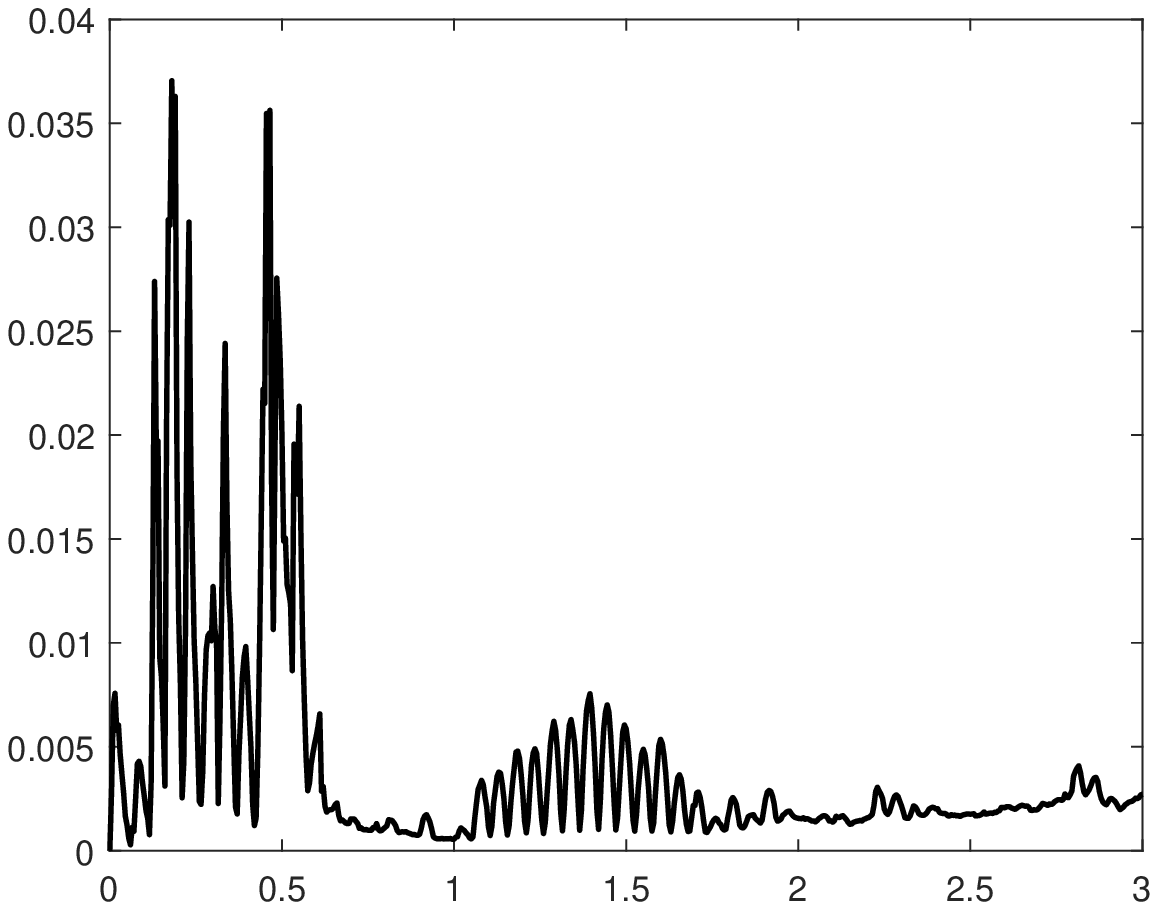}}
\caption{Stochastic generalized  Kuramoto--Sivashinsky (SgKS) equation (\ref{eq:1}): (a) absolute error $\Delta_a u$ as a function of time $t$, and (b) relative error $\Delta_r u$ as a function of time $t$.}
\label{fig:16}
\end{figure}

Some results are presented in Figures \ref{fig:14}-\ref{fig:16}. These figures show that the semi-analytical solution and the WCE approximate numerical solution are in good agreement, and their absolute difference $\Delta_a u$ is, in general, order $10^{-3}$ on the time interval $[0,3]$. Figure \ref{fig:16}(a) shows that the absolute error $\Delta_a u\sim O(\varsigma t)$ as in test problems 1,2 and 3, and where $\varsigma \sim O(10^{-3})$. It is seen in Figure \ref{fig:16}(b) that the relative (error) $\Delta_r u$ is $O(10^{-2})$ or less  for all $t$ in the interval $[0,3]$ and its maximum value on this interval is $4\%$. 

\section{Discussions and concluding remarks}\label{sec:6}
We have computed and examined the WCE based approximate numerical solutions to the stochastic Kuramoto-Sivashinsky (SKS) equation and stochastic generalized Kuramoto-Sivashinsky (SgKS) equation with Brownian motion forcing. A semi-analytical solution procedure was discussed as well. Some preliminary tests involving the linear (SgKS) equation were performed.

In our preliminary tests in section \ref{sec:3}, we linearized the stochastic generalized Kuramoto-Sivashinsky equation and solved the linear stochastic evolution equation using the WCE based numerical method. We considered some IBVPs for which analytical solutions can be derived and and we have performed four realizations. In each realization, the WCE approximate numerical solution was contrasted with the analytical solution, and the accuracy of the results was examined following our analytical predictions in section \ref{subsec:3.2}. It was found that there is a good agreement between the WCE based numerical solutions and the analytical solutions. The absolute and relative errors increased with time and the relative error attained a maximum value of 1\% over the time interval $[0,3]$. The errors can be minimized by using a small time step size and by increasing the order of the WCE as predicted by formula (\ref{eq:31}).

Next, some IBVPs involving the SKS and SgKS equation were considered. However, in this case, there are no analytical solutions to compare with. WCE based numerical solutions were compared with semi-analytical solutions obtained using the procedure described in section \ref{subsec:4.2}. Non-homogeneous stochastic Dirichlet boundary conditions were implemented at the boundaries of the domain to take into consideration the stochastic evolution of the solutions near the boundaries, and four test problems were considered. 

In each test problem, the WCE based numerical solution was contrasted with the semi-analytical solution, and it was found that there is good agreement between the WCE based numerical solutions and semi-analytical solutions. The absolute error (difference) and the relative difference were evaluated using formulas (\ref{eq:30}) and (\ref{eq:32}). It was found that the absolute error is $O(\varsigma t)$, where $\varsigma$ is a small constant of order $10^{-3}$, and the relative error is order $10^{-2}$ over time interval $[0.3]$.

The results presented in section \ref{sec:5} were obtained using the predictor-corrector method described in section \ref{subsec:4.1}, a time step size $\Delta t=0.005$ and a spatial step size $\Delta x=0.2$. The predictor-corrector method is effective and quite simple and was, for example, used in \cite{Nijimbere2014ionospheric,Nijimbere2016nonlinear} to solve the vorticity equations. Additional computations were done using the time step sizes $\Delta t=0.01,0.02, 0.05$ and different $\Delta x=0.02, 0.05, 0.1$. It was found that the accuracy is higher when the ratio $\Delta t/\Delta x\ll1$. 

It is worth to point out that our solutions are random fields which are characterized by their means and variances. Once it has been verified that the WCE based numerical methods are effective, it is then straight forward to compute the mean and the variance of the solution using the formulas in (\ref{eq:15}).

In conclusion, the results of this study illustrate that the WCE based numerical methods are powerful methods for solving stochastic evolution partial differential equations (PDEs) such as the stochastic Kuramoto-Sivashinsky equation driven by Brownian motion forcing.

%\bibliographystyle{IEEEtran}
%
%\bibliography{paper}

\end{document}